\documentclass[12pt]{amsart} 

\usepackage[margin=1in]{geometry} 

\usepackage{amssymb,amsmath,amsthm,mathrsfs,enumerate,multicol,esint,color}

\allowdisplaybreaks

\title[Weighted embeddings for Hermite spaces]{Weighted embeddings for function spaces associated with  Hermite expansions}

\author{The Anh Bui}
\address{Department of Mathematics, 
Macquarie University, 
NSW 2109, Australia}
\email{the.bui@mq.edu.au}

\author{Ji Li}
\address{Department of Mathematics, Macquarie University, 
NSW 2109, Australia}
\email{ji.li@mq.edu.au}

\author{Fu Ken Ly}
\address{School of Mathematics and Statistics, Faculty of Science \& The Mathematics Learning Centre, Education Portfolio, The University of Sydney, 
NSW 2006, Australia}
\email{ken.ly@sydney.edu.au}

\thanks{2010 {\it Mathematics Subject Classification}: 46E35, 42C15, 33C45. }
\thanks{{\it Key words and phrases}: Embedding, Hermite expansions, Frames, Weights}
\thanks{\copyright \;2021. This manuscript version is made available under the CC-BY-NC-ND 4.0 license {http://creativecommons.org/licenses/by-nc-nd/4.0/}}

\newcommand{\RR}{\mathbb{R}} 
\newcommand{\NN}{\mathbb{N}} 
\newcommand{\ZZ}{\mathbb{Z}} 
\newcommand{\sz}{\mathscr{S}} 
\newcommand{\MM}{\mathcal{M}}

\newcommand{\f}{\frac}
\newcommand{\lesi}{\lesssim}

\newcommand{\LL}{\mathcal{L}}
\newcommand{\PP}{\mathbb{P}} 
\newcommand{\QQ}{\mathbb{Q}} 
\newcommand{\ip}[1]{\langle #1 \rangle} 
\newcommand{\E}{\mathcal{E}} 
\newcommand{\X}{\mathcal{X}} 
\newcommand{\Ind}{\mathbf{1}} 
\newcommand{\Q}{\mathcal{Q}}
\newcommand{\A}{\mathcal{A}}
\newcommand{\B}{\mathcal{B}}
\newcommand{\AH}{\mathbb{A}} 
\newcommand{\wt}[1]{\widetilde{#1}} 
\newcommand{\ve}{\varepsilon}
\newcommand{\diff}{\triangle} 
\newcommand{\floor}[1]{\lfloor #1 \rfloor}
\newcommand{\ceil}[1]{\lceil #1 \rceil}
\newcommand{\ep}{\epsilon}

\newcommand{\vph}{\varphi}
\newcommand{\embed}{\hookrightarrow}
\newcommand{\e}{\mathrm{e}}
\newcommand{\I}{\mathcal{I}}
\newcommand{\cro}{\varrho}

\DeclareMathOperator{\supp}{supp\,} 

\theoremstyle{plain}
\newtheorem{Theorem}{Theorem}[section]
\newtheorem{Lemma}[Theorem]{Lemma}
\newtheorem{Proposition}[Theorem]{Proposition}

\newtheorem{Definition}[Theorem]{Definition}
\theoremstyle{definition}
\newtheorem{Remark}[Theorem]{Remark}

\theoremstyle{remark}

\numberwithin{equation}{section}	

\def\barint{\kern4pt
\raise3.4pt\hbox{\vrule height.8pt width5pt}%
\kern-9pt 
\int}

\newcommand{\aver}{-\hskip-0.46cm\int}

\def\XXint#1#2#3{{\setbox0=\hbox{$#1{#2#3}{\int}$}
     \vcenter{\hbox{$#2#3$}}\kern-.5\wd0}}

\begin{document}

\begin{abstract}
We study weighted Besov and Triebel--Lizorkin spaces associated with Hermite expansions and obtain (i) frame decompositions, and (ii) characterizations of continuous Sobolev-type embeddings. The weights we consider generalize the Muckhenhoupt weights.
\end{abstract}

\maketitle 
\tableofcontents

\section{Introduction}
One of the most important features of smooth function spaces is the presence of Sobolev type embeddings. The ability to exchange regularity for increased integrability is a critical tool in the study of partial differential equations. In particular weighted embedding theorems for Besov and Triebel--Lizorkin spaces have found applications in a variety of situations (see \cite{MV1,MV3} and the references therein).

To facilitate the rest of this discussion let us recall such embeddings. In what follows, by a weight $w$ we shall mean that $w$ is a locally integrable non-negative function defined  on $\mathbb{R}^n$. For a measurable subset $E\subset \mathbb R^n$ and a weight $w$, we write
$$w(E)=\int_E w(x)dx.$$ 
Given a weight $w$  let  $B^{p,q}_{\alpha,w}$ and $F^{p,q}_{\alpha,w}$ be the weighted inhomogeneous Besov and Triebel--Lizorkin spaces on $\RR^n$ respectively, where $0<p<\infty$ is the integrability index, $0<q<\infty$ the fine index, and $\alpha\in\RR$ the smoothness index. Suppose $\alpha_2\ge\alpha_1$, $0<p_2\le p_1<\infty$, $0<q_1,q_2\le \infty$ and $\alpha_1-\gamma/p_1 = \alpha_2-\gamma/p_2$ for some $\gamma>0$. Then under suitable conditions on $w$ one has the continuous embedding
\begin{align}\label{eq:Femb}
F^{p_2,q_2}_{\alpha_2,w}\embed F^{p_1,q_1}_{\alpha_1,w},
\end{align}
and for $0<q_2\le q_1\le \infty$,
\begin{align}\label{eq:Bemb}
B^{p_2,q_2}_{\alpha_2,w}\embed B^{p_1,q_1}_{\alpha_1,w}.
\end{align}
Situations for which \eqref{eq:Femb}--\eqref{eq:Bemb} hold true include of course the unweighted case of $w=1$ and $\gamma=n$, which can be found in \cite{Franke,Jawerth,Triebel}. Other situations include power weights \cite{HS,MV1} and radial weights \cite{DDS}.  When $w$ is a Muckenhoupt weight, \eqref{eq:Femb}--\eqref{eq:Bemb} are known to hold if and only if $w$ satisfies a kind of (local) lower bound property, namely that
\begin{align}\label{eq:LLB}
w(B(x,r))\ge C r^\gamma,\qquad x\in\RR^n,\; 0<r\le 1
\end{align}
for some $C>0$ independent of $x$ and $r$, where $B(x,r)=\{y\in \mathbb R^n: |x-y|<r\}$.
This was proved in \cite[Theorem 2.6]{Bui}, \cite[Proposition 2.1(i)]{HS} and \cite[Theorem 1.2]{MV2}. That latter two articles in fact consider two-weight embeddings of which \eqref{eq:Femb}--\eqref{eq:Bemb} are a special case.

In recent times various works have been devoted to generalizations of the Besov and Triebel--Lizorkin spaces in a variety of directions. Such directions include generalizing the ambient space $\RR^n$ \cite{DH,HanSawyer}, replacing the Fourier system by other orthogonal expansions \cite{Dz,Ep,KPPX,KPX1,PX}, adapting the function spaces to  a general class of operators \cite{BBD,BDY,KP,GKKP} or hybrids of the above. 

It is thus natural to ask if embeddings similar to \eqref{eq:Femb}--\eqref{eq:Bemb} can be retained in these generalizations. 
For the setting of a space of homogeneous type, a complete characterization of \eqref{eq:Femb}--\eqref{eq:Bemb} was given in \cite{HHL} by \eqref{eq:LLB}, where $w=\mu$ is the underlying measure, unifying and completing several earlier results \cite{HYC1,HYC2,HYS1,HL}.

\bigskip

In the present article we consider the case of Hermite function expansions \cite{Th}, which are naturally associated with the harmonic oscillator
$$ \LL=-\Delta+|x|^2 $$
on $\RR^n$ for $n\ge 1$. Besov and Triebel--Lizorkin spaces associated with Hermite expansions have been developed in \cite{BD,Dz,Ep,PX}.
In this article we introduce weighted versions with the aim of characterizing Sobolev--type embeddings. 

To describe our main results let us first define our function spaces. 
Let $(\vph_\circ, \varphi)$ be an admissible pair as defined in Definition \ref{def: admissible}. Let $\{h_\xi\}_{\xi\in \NN_0^n}$ be the multi-dimensional Hermite functions (see Section \ref{sec:prelim}) and consider, for $j\ge 0$, the Littlewood--Paley type operators defined by
$$\vph_j(\sqrt{\LL})  f := \sum_{\xi\in\NN_0^n} \vph_j(\sqrt{2|\xi|+n})\ip{h_\xi,f} h_\xi,\qquad f\in\mathscr{S}'$$
where $\ip{h_\xi,f}:=f(h_\xi)$. Then we introduce the following weighted spaces. 
Let  $\alpha\in\RR$, $0<q\le\infty$  and $w$ be a weight.
For $0<p\le\infty,$ we define the  weighted Hermite Besov space $B^{p,q}_{\alpha,w}$ as the class of tempered distributions $f\in \sz'(\RR^n)$ such that 
\begin{align}\label{eq:besov} 
\Vert f\Vert_{B^{p,q}_{\alpha,w}} = \Big(\sum_{j\in\NN_0}\big(2^{j\alpha}\Vert \varphi_j(\sqrt{\LL})f\Vert_{L^p_w}\big)^q\Big)^{1/q}<\infty,
\end{align}
while for $0<p<\infty,$  the weighted Hermite Triebel--Lizorkin space  $F^{p,q}_{\alpha,w}$ consists of all those $f\in \sz'(\RR^n)$ for which
\begin{align}\label{eq:triebel}
 \Vert f\Vert_{F^{p,q}_{\alpha,w}} = \Big\Vert \Big(\sum_{j\in\NN_0}\big(2^{j\alpha}| \varphi_j(\sqrt{\LL})f|\big)^q\Big)^{1/q} \Big\Vert_{L^p_w}<\infty.
 \end{align}
We adopt the usual modifications to the norm when $q=\infty$ (or $p=\infty$ in the Besov case). The spaces can be shown to be independent of our choice of $(\vph_\circ,\vph)$ (see Theorem \ref{thm:frames}).

Our aim in this paper is to characterize embedding theorems for such spaces for a suitable class of weights (which we denote by $\AH_\infty(\LL)$). These weights  are more general than (and include) the Muckenhoupt weights. 
Such weights, introduced in \cite{BHS}, have sparked some recent activities (see for example  \cite{BCH,LRW,Ly,Tang}). They can be non-doubling and may even admit a certain degree of exponential growth. Note that embedding in the unweighted situation (with $w=1$) was given in \cite{PX}. More details about the kinds of weights considered in this paper can be found in Section \ref{sec:weights}. 

Our main result, Theorem \ref{thm:embedding1} characterizes \eqref{eq:Femb}--\eqref{eq:Bemb} via a similar condition to \eqref{eq:LLB}.
For the definition of the weight class  $\AH_\infty(\LL)$ we refer to  Definition \ref{def:weights}.
\begin{Theorem}\label{thm:embedding1}
Let $\gamma>0$ and $w\in \AH_\infty(\LL)$. Then one has the embeddings \eqref{eq:Femb} or \eqref{eq:Bemb}  if and only if $w$ satisfies the following ``Hermite lower bound property'' of order $\gamma$:
there exists a constant $C>0$ such that,
\begin{align}\label{eq:HLB}
w(B(x,r))\ge C\,r^\gamma,\qquad  x\in\RR^n,\quad 0<r\le \f{1}{1+|x|_\infty}.
\end{align}
\end{Theorem}
Examples of weights which satisfy \eqref{eq:HLB} include the power weights $w(x)=|x|^\ep$, for $-n<\ep<\infty$ provided that $\gamma\ge n+\ep$. An example that is not a Muckenhoupt weight is the exponential weight $w(x)=e^{\ep|x|^2}$ with $\gamma\ge n$.

Some remarks on our method of proof are in order. Our overall approach to Theorem \ref{thm:embedding1} is to obtain embeddings for the associated sequence spaces (see Definition \ref{def:sequences} and Theorem \ref{thm:embedding2}) and then pass to the function spaces using suitable frame characterizations (Theorem \ref{thm:frames}). To obtain embeddings in the Triebel--Lizorkin scale we employ the well known distribution function method which harks back to Jawerth \cite{Jawerth} and refined in \cite{Triebel}. The embeddings on the Besov scale are relatively easier. 

We give a final important remark on the relevance of condition \eqref{eq:HLB}. Our proofs rely on the ability to pass from a lower bound estimate on frame elements to the estimate on balls (Proposition \ref{prop: lb}). This is a simple matter in the classical situation since the frame elements are dyadic cubes  whose sizes depend only on the scale. In the Hermite setting however, the frame elements are so-called `tiles' (see Section \ref{sec:tiles}) which are dependent on both scale \emph{and} location. This principle is quantified in \eqref{eq:HLB} and Proposition \ref{prop: lb}.

The organization of the article is as follows. In Section \ref{sec:prelim} we give the necessary background that will be used throughout the paper including facts concerning Hermite functions, kernel estimates for our Littlewood-Paley type operators, Hermite weights and Hermite tiles (frame elements).
 In particular the associated weighted sequence spaces are introduced here in Definition \ref{def:sequences}. Suitable frames  are introduced in \ref{sec:frames}   and a frame characterization (Theorem~\ref{thm:frames}) is proved. Our embedding results are proved in Section \ref{sec:embedding}. Here we first develop embeddings for Hermite sequence spaces (Theorem \ref{thm:embedding2}) before finally giving the proof of Theorem \ref{thm:embedding1}.

As a final remark, it would be interesting to extend our  results to two-weights embeddings, as done in \cite{HS,MV2}. However we do not pursue this here.

{\bf Notation:} We conclude this introduction with some notational matters. 
We set $|x|_\infty:=\max\limits_{1\le i\le n} |x_i|$.
We denote by $Q(x,r):=\{ y\in\RR^n: |x-y|_\infty<r\}$ the cube centered at $x$ with sidelength $2r$. By a `cube $Q$' we mean the cube $Q(x_Q,r_Q)$ with some fixed centre $x_Q$ and sidelength $2r_Q$.
We denote $\lambda_k:=2k+n$ for $k\in\NN_0$, where $\NN_0=\NN\cup \{0\}$.
For two Banach spaces $X_0$ and $X_1$ the notation $X_0\embed X_1$ means that $X_0$ embeds continuously into $X_1$. We denote by $\sz(\RR^n)$ the space of Schwartz functions on $\RR^n$ and by $\sz'(\RR^n)$ the space of tempered distributions on $\RR^n$. The letter $n$ will always mean the Euclidean dimension.

\section{Preliminaries}\label{sec:prelim}
In this section we give the necessary  material that will be used throughout the rest of the paper.   After detailing some basic facts concerning Hermite functions we introduce and give some kernel estimates for our Littlewood--Paley operators in Section \ref{sec:LP}. In Section \ref{sec:weights} we define and give some useful estimates for the Hermite weight classes employed throughout the paper. Finally in Section \ref{sec:tiles} we detail the construction of dyadic `tiles' and introduce weighted sequence spaces associated with Hermite expansions.

For each $k\in\NN$, the Hermite function of degree $k$ is 
\begin{align*}
	h_k(t)=(2^k k! \sqrt{\pi})^{-1/2}H_k(t)e^{-t^2/2}, \qquad t\in\RR{},
\end{align*}
where 
	$H_k(t)=(-1)^ke^{t^2}\partial_t^k(e^{-t^2})$
is the $k$th Hermite polynomial. 

The $n$-dimensional Hermite functions $h_\xi$ are defined over the multi-indices $\xi$ by 
\begin{align*}
	h_\xi(x)=\prod_{j=1}^n h_{\xi_j}(x_j),\qquad x\in\RR^n, \quad \xi \in \NN^n_0.
\end{align*}
The Hermite functions are eigenfunctions of $\LL$ in the sense that $\LL(h_\xi) = \lambda_{|\xi|} h_\xi, $
where $\lambda_k=2k+n$. Furthermore they form an orthonormal basis for $L^2(\RR^n)$. 

We also set, for each $j=1,\ldots, n$,
\[
A_j =-\frac{\partial}{\partial x_j} +x_j \qquad\text{and}\qquad A^*_j=\frac{\partial}{\partial x_j} +x_j,
\]
to be the so-called `creation' and `annihiliation' operators respectively.
Then it follows that Harmonic oscillator $\LL=-\Delta+|x|^2$ satisfies
\[
\mathcal L =\frac{1}{2}\sum_{j=1}^n\Big[A_j^*A_j+ A_jA^*_j\Big].
\]
Let $W_k=\text{span}\{h_\xi: |\xi|=k\}$ and $ V_N = \bigoplus^N_{k=0} W_k$. We define the orthogonal projection of $f$ onto $W_k$ by
\begin{align*}
	\PP_kf = \sum_{|\xi|=k} \ip{f, h_\xi} h_\xi 
	\qquad\text{with kernel}\qquad \PP_k(x,y)=\sum_{|\xi|=k} h_\xi(x)h_\xi(y).
\end{align*}
We also define the orthogonal projection of $f$ onto $V_N$ by
\begin{align*}
	\QQ_N f = \sum_{k=0}^N \PP_kf
	\qquad\text{with kernel}\qquad
	\QQ_N(x,y)=\sum_{k=0}^N \sum_{|\xi|=k}h_\xi(x)h_\xi(y).
\end{align*}
The following bounds are known (see \cite[p. 376]{PX}): there exists $\vartheta>0$ such that for any $N\ge 1$
\begin{align}\label{QQ bound}
	\QQ_{N}(x,x) 
	\lesssim \left
\lbrace 
	\begin{array}{ll}
			N^{n/2}  \qquad &\forall x\\
			e^{-2\vartheta|x|^2} \qquad &\text{if}\quad |x| \ge \sqrt{4N+2}.
				\end{array}
\right.
\end{align}
We also set
\begin{equation}\label{eq-eN}
\begin{aligned}
	\e_{N}(x) := \left
\lbrace 
	\begin{array}{ll}
			1  \qquad &\text{if}\quad |x|^2  < N\\
			e^{-\vartheta |x|^2} \qquad &\text{if}\quad |x|^2\ge N.
	\end{array}
\right.
\end{aligned}
\end{equation}
Then it follows  from \eqref{QQ bound} that for any $\ve>4$ and $N\in\NN$, there exists $C>0$ depending only on  $N,$ $n,$ $\ve$ and $\vartheta$ such that 
\begin{align}\label{QQ est}
\QQ_{4^j+N}(x,x)\le C\; 2^{jn}(\e_{\ve 4^j}(x))^2 \qquad \forall \; j\in\NN_0.
\end{align}

\subsection{Littlewood--Paley operators}\label{sec:LP}
The following band-limited functions will play a fundamental role in this paper. Recall that they were used to define our weighted function spaces in \eqref{eq:besov} and \eqref{eq:triebel}.
\begin{Definition}[Admissible functions]\label{def: admissible}
We say that $(\vph_\circ, \varphi)$ are an admissible pair if $\vph_\circ, \varphi\in C^\infty(\RR_+)$ and
\begin{align*}
	\supp\varphi_{\circ}\subset [0,2], && |\varphi_\circ| > c> 0 \quad\text{on}\quad [0,2^{-3/4}], && \varphi_\circ^{(m)}(0)=0, \quad \forall m\in\ZZ_+ \\
	\supp\varphi\subset [\tfrac{1}{2},2], && |\varphi|>c>0 \quad\text{on}\quad [2^{-3/4},2^{3/4}].
\end{align*}
Given an admissible pair $(\varphi_\circ, \varphi)$, we set $\varphi_j(\lambda):=\varphi(2^{1-j}\lambda)$ if $j\ge 1$ and $\vph_0(\lambda)=\vph_\circ(2\lambda)$  and call the resulting collection $\{\varphi_j\}_{j\ge 0}$ an admissible system. 
\end{Definition}
Since the Hermite functions $\{h_\xi\}_{\xi\in\NN_0^n}$ are members of $\mathscr{S}(\RR^n)$, then for an admissible system $\{\vph_j\}_{j\in\NN_0}$ we may define the operators $\vph_j(\sqrt{\LL})$ on $\mathscr{S}'(\RR^n)$ by 
\begin{align*}
	\vph_j(\sqrt{\LL})  f = \sum_{\xi\in\NN_0} \vph_j(\sqrt{\lambda_\xi})\,\PP_\xi f, \qquad f\in \mathscr{S}'(\RR^n),
\end{align*}
where $\ip{f,\phi}=f(\phi)$ for $f\in \mathscr{S}'(\RR^n)$. The kernels of operators $\vph_j(\sqrt{\LL})$ are given by 
\begin{align*}
	\vph_j(\sqrt{\LL})(x,y)=\sum_{\xi\in\NN_0} \vph_j(\sqrt{\lambda_\xi})\PP_\xi(x,y) = \sum_{\xi\in\NN_0}\vph_j(\sqrt{\lambda_\xi})\sum_{|\mu|=\xi}h_\mu(x)h_\mu(y).
\end{align*}
Let $\{\I_j\}_{j\in\NN_0}$ denote the following subsets of $\NN_0$: $\I_j=\big[\f{1}{2}4^{j-2}-\floor{\frac{n}{2}}, \f{1}{2}4^{j}-\ceil{\frac{n}{2}}\big]\cap \NN_0$ for $j\in \NN,$ $\I_0=\{0\}$ if $n=1$  and $\I_0=\emptyset$ if $n\ge 2$. 
The support of $\varphi_j$ implies that
\begin{align}\label{eq:kernel}\varphi_j(\sqrt{\LL})(x,y)=\sum_{\xi\in \I_j}\varphi_j(\sqrt{\lambda_\xi})\, \sum_{|\mu|=\xi}h_\mu(x)h_\mu(y).\end{align}

We give now some estimates on the kernels of the operators $\vph_j(\sqrt{\LL})$. Some of these have previously appeared in \cite{LN,PX} but we give their proofs for the convenience of the reader.

\begin{Proposition}[Kernel estimates]\label{prop: kernel}
Let $\{\varphi_j\}_{j\ge 0}$ be an admissible system and fix $N\ge 1$ and $\ve>4$. 
Then there exists $C>0$ depending on $N$, $\vartheta$, $\ve$, and  {\upshape($\vph_\circ$, $\vph$) }such that
\begin{align}\label{eq:kernelbound}
	|\varphi_j(\sqrt{\LL})(x,y)|\le C \f{2^{jn}}{(1+2^j|x-y|)^N} \e_{\ve 4^j}(x) \e_{\ve 4^j}(y),\qquad \forall\;x,y\in\RR^n.
\end{align}
Suppose that $\{\psi_j\}_{j\ge 0}$ is another admissible system. 
	If $|j-k|\ge 3$ then
	\begin{align}\label{eq:orthog1}
	\varphi_j(\sqrt{\LL})  \psi_k(\sqrt{\LL})(x,y) = 0, \qquad \forall\;x,y\in\RR^n.
	\end{align}
	There exists $C>0$ depending on $N$, $\vartheta$, $\ve$, {\upshape($\vph_\circ$, $\vph$) } and {\upshape($\psi_\circ$,$\psi$)} such that whenever $|j-k|\le 3$ then 
	\begin{align}\label{eq:orthog2}
	| \varphi_j(\sqrt{\LL})  \psi_k(\sqrt{\LL})(x,y) | \le C \f{2^{kn}}{(1+2^k|x-y|)^N}\e_{\ve 4^{j}}(x) \e_{\ve 4^{k}}(y)
	\end{align}
for every $x,y\in\RR^n$. 
\end{Proposition}
\begin{proof}
We begin with the proof of \eqref{eq:kernelbound}.
 If $|x-y| \le 2^{-j}$ then by the Cauchy--Schwarz inequality,
\begin{align*}
|\varphi_j(\sqrt{\LL})(x,y)|
\le \|\vph_j\|_\infty \sum_{\xi\in\I_j}  \sum_{|\mu|=\xi}|h_{\mu}(x)|\,|h_\mu(y)| 
\lesi  \QQ_{4^j+N}(x,x)^{1/2} \QQ_{4^j+N}(y,y)^{1/2}.
\end{align*}
If on the other hand, $|x-y| > 2^{-j}$,
then applying the identity in \cite[Lemma 8]{PX} (see also \cite[p. 72]{Th} or \cite[(B.1)]{LN}) to \eqref{eq:kernel} we obtain
\begin{align}\label{eq:kernel identity}
\big|2^N (x_i-y_i)^N \varphi_j(\sqrt{\LL})(x,y)\big|
	&= \sum_{N/2\le \ell\le N} c_{\ell,N} \sum_{\xi\in\I_j} \big|\diff^\ell \varphi_j(\sqrt{\lambda_\xi})\big|\,\big|(A_i^{(y)}-A_i^{(x)})^{2\ell-N} \PP_\xi(x,y)\big|,
\end{align}
where $c_{\ell,N} = (-4)^{N-\ell}(2N-2\ell-1)!!\binom{N}{2\ell-N}$ and $A_i^{(x)}$ and $A_i^{(y)}$ denote the operator $A_i$ applied to the $x$ and $y$ variables, respectively. Here the symbol $\diff$ denotes the forward difference operator given by $\diff f(\xi)=f(\xi+1)-f(\xi)$.

Using the mean value theorem for finite differences along with Hoppe's chain rule formula we have
\begin{align}\label{eq:hoppe}
	|\diff^\ell \varphi_j(\sqrt{\lambda_\xi})| 
	= |\partial^\ell_\nu \varphi_j(\sqrt{\lambda_\nu})|  
	= \Big| \sum_{r=1}^\ell c_r\,\varphi_j^{(r)}(\sqrt{\lambda_\nu}) \lambda_\nu^{r/2-\ell}2^{-jr}\Big|, \qquad |c_r|\le \ell!
\end{align}
for some $\xi<\nu<\xi+\ell$. Now our assumptions on $\vph_\circ$ and $\vph$ imply $\varphi_j^{(m)}(0)=0$ for every $m\in \ZZ_+$ and, in conjunction  with Taylor's remainder theorem,
gives 
$$|\varphi_j^{(r)}(x)|\le \max\big\{\Vert \vph_\circ^{(N)} \Vert_{\infty}, \Vert \vph^{(N)}\Vert_\infty\big\}  \;|2^{1-j}x|^{N-r}$$
 for any $N\ge 1$ and $x\in \RR_+$. Inserting this into \eqref{eq:hoppe} we arrive at
\begin{align}\label{taylor}
	|\diff^\ell \varphi_j(\sqrt{\lambda_\xi})| 
\lesssim  \lambda_\xi^{N/2-\ell}2^{-jN} 
	\end{align}
since $\xi<\nu<\xi+\ell$, with an implicit constant that depends on $\vph_\circ$, $\vph$, $\ell$ and $N$.
Next we recall the following estimate from \cite[pp. 397-398]{PX},
\begin{align}\label{PX est}
|(A_i^{(y)}-A_i^{(x)})^{2\ell-N} \PP_\xi(x,y)| 
\lesssim \lambda_{\xi}^{\ell-N/2} \Big(\sum_{k=0}^{2\ell-N}\PP_{\xi+k}(x,x)\Big)^{1/2} \Big(\sum_{k=0}^{2\ell-N}\PP_{\xi+k}(y,y)\Big)^{1/2},
\end{align}
which is obtained through the Binomial theorem and Cauchy--Schwarz inequality. 
Inserting \eqref{taylor} and \eqref{PX est} into \eqref{eq:kernel identity} and applying Cauchy--Schwarz inequality we obtain
\begin{align*}
 & |(x_i-y_i)^N \varphi_j(\sqrt{\LL})(x,y)| \\
 &\qquad\qquad\lesssim 2^{-jN} \sum_{N/2\le \ell\le N}  c_{\ell,N}\sum_{\xi\in\I_j}  \Big(\sum_{k=0}^{2\ell-N}\PP_{\xi+k}(x,x)\Big)^{1/2} \Big(\sum_{k=0}^{2\ell-N}\PP_{\xi+k}(y,y)\Big)^{1/2} \\
 &\qquad\qquad\lesssim  2^{-jN} \QQ_{4^j+N}(x,x)^{1/2} \QQ_{4^j+N}(y,y)^{1/2}.
\end{align*}

Now combining both cases we arrive at
$$ |\varphi_j(\sqrt{\LL})(x,y)|\lesssim \f{\QQ_{4^j+N}(x,x)^{1/2} \QQ_{4^j+N}(y,y)^{1/2} }{(1+2^j|x-y|)^N},$$
which, in view of \eqref{QQ est}, yields \eqref{eq:kernelbound}.

We turn to the proof of \eqref{eq:orthog1}.
By orthonormality of the Hermite functions
\begin{align*}
	\varphi_j(\sqrt{\LL}) \psi_k(\sqrt{\LL})(x,y)
	&=\sum_{\xi\in\NN_0}\varphi_j(\sqrt{\lambda_\xi})\varphi_k(\sqrt{\lambda_\xi}) \sum_{|\mu|=\xi}h_\mu(x) h_\mu(y).
\end{align*}
From the supports of $\vph_j$ and $\psi_k$, it can be easily seen that
$\varphi_j(\sqrt{\lambda_\xi})\psi_k(\sqrt{\lambda_\xi})=0$ whenever $|j-k|\ge 3$, which gives \eqref{eq:orthog1}.

\medskip

We now prove \eqref{eq:orthog2}. 
From the bound \eqref{eq:kernelbound} and the fact that $|j-k|\le 3$ we have
\begin{align*}
| \varphi_j(\sqrt{\LL})  &\psi_k(\sqrt{\LL})(x,y) | \\
&\le  \int \big|\varphi_j(\sqrt{\LL})(x,z)\big| \,\big|\psi_k(\sqrt{\LL})(z,y)\big|\,dz  \\
&\lesssim \e_{\ve 4^j}(x) \e_{\ve 4^k}(y)\int \f{2^{kn}}{(1+2^k|x-z|)^N}\f{2^{kn}}{(1+2^k|z-y|)^M} \,dz 
\end{align*}
for some $M>N+n$.
The triangle inequality then gives
\begin{align*}
| \varphi_j(\sqrt{\LL})  &\psi_k(\sqrt{\LL})(x,y) | 
\lesi \f{2^{kn}\e_{\ve 4^j}(x) \e_{\ve 4^k}(y)}{(1+2^k|x-y|)^N} \int \f{2^{kn}}{(1+2^k|z-y|)^{M-N}} \,dz.
\end{align*}
Finally we observe that since $M>N+n$ then the integral is finite which completes our proof.
\end{proof}

\subsection{Hermite weights}\label{sec:weights}
Here we define the classes of weights used in our function and sequence spaces. They generalize the Muckenhoupt classes and were introduced in \cite{BHS}.
These weights are based on  the following `critical radius' function.
\begin{align}\label{rho}
	\cro(x):=\f{1}{1+|x|_\infty},
\end{align}
which was introduced in \cite{Shen} in a more general context. It can be easily seen that there exists $C>0$ and $\kappa\ge 1$ such that
\begin{align}\label{slow growth}
	C^{-1}\cro(x)\Big(1+\f{|x-y|}{\cro(x)}\Big)^{-\kappa}\le \cro(y)\le C\cro(x)\Big(1+\f{|x-y|}{\cro(x)}\Big)^{\f{\kappa}{\kappa+1}}
\end{align}
for any $x,y\in\RR^n$ (see also \cite[Lemma 1.4]{Shen}).
\begin{Definition}[Hermite weights]\label{def:weights}
For each cube $Q$ and $\theta\ge0$ set
\begin{align*}
	\Psi_\theta(Q):= \Big(1+\f{\ell(Q)}{\cro(x_Q)}\Big)^\theta,
\end{align*}
where $\ell(Q)$ is the side-length of the cube $Q$.

For each $1<p<\infty$ and $\eta\ge0$ we define the class $\AH_p^\eta(\LL)$ to be the collection of all non-negative and locally integrable functions $w$ such that for some $C>0$,
\begin{align}\label{weights1}
	\Big(\aver_Q w(y)\,dy\Big)^{1/p}\Big(\aver_Q w^{1-p'}(y)\,dy\Big)^{1/p'}\le C\,\Psi_\eta(Q)
\end{align}
for every cube $Q$. 
For $p=1$ we require the following in place of \eqref{weights1}
\begin{align*}
	\f{1}{\Psi_\eta(Q)}\aver_Q w(y)\,dy\le w(x)\qquad \text{a.e.}\quad x\in Q.
\end{align*}
We also set
$$ \AH_p(\LL)=\bigcup_{\eta\ge0}\AH_p^\eta(\LL) \qquad\text{and}\qquad  \AH_\infty(\LL) :=\bigcup_{ p\ge 1}\AH_p(\LL).$$
\end{Definition}

\begin{Remark}\label{rem: weights}
\begin{enumerate}[(i)]
	\item Throughout the rest of the article we shall drop the $\LL$ in the notation and just write $\AH_p^\eta$ in place of $\AH_p^\eta(\LL)$, and $\AH_p$ in place of $\AH_p(\LL)$. 
	\item We denote the usual class of Muckenhoupt weights by $A_p$ for $p\ge 1$. Then one may observe that the Hermite weights form a larger class than the Muckenhoupt weights in the sense that for each $p\ge 1$ and $\eta\ge 0$ we have $A_p\subset \AH_p^\eta$. 
	\item We note that these weights are increasing in both parameters in the sense that for any $\eta\ge 0$, then whenever $1\le p_1\le p_2<\infty$ we have
	$$ \AH_1^\eta \subset \AH_{p_1}^\eta\subset \AH_{p_2}^\eta $$
	and for any $p\ge 1$ and $\eta_1\le \eta_2$ then
	$$ \AH_p^{\eta_1}\subset \AH_p^{\eta_2}.$$
	\item Examples of $\AH_\infty$ weights are the following. For each $p>1$ we have $e^{\epsilon_1|x|^2}, e^{-\epsilon_2|x|^2}\in\AH_p$ for some $\epsilon_1, \epsilon_2>0$. See \cite{LRW}.
	\item The class $\AH_\infty$ is `open' in the sense that if $w\in\AH_p$, then $w\in\AH_{p-\epsilon}$ for some $\epsilon>0$. Thus for each $w\in \AH_\infty$ we may define
	\begin{align*}
		r_w:=\inf\big\{r>1:\; w\in\AH_r\big\}.
	\end{align*}
	\end{enumerate}
\end{Remark}

We will also need to employ certain maximal functions throughout our paper. 
For each $s>0$ and $\theta\ge 0$ we define a maximal function as follows:
\begin{align}\label{mf}
	\MM^{\theta}_s f(x):=\sup_{Q\ni x} \Big(\f{1}{\Psi_\theta(Q)}\aver_Q |f(y)|^s\,dy\Big)^{1/s}.
\end{align}
When $s=1$ we drop the $s$ and write $\MM^\theta$. 
These maximal functions satisfy a (weighted) Fefferman--Stein inequality for vector valued functions.
\begin{Lemma}[Fefferman--Stein type inequality]\label{lem:FS} 
Let $0<p,q<\infty$ and $\{f_j\}_{j\in \NN_0}$ be a sequence of locally integrable functions defined on $\RR^n$. 
\begin{enumerate}[\upshape(a)]
\item If $\min\{p,q\}>1$ and  $w\in\AH_p^\eta$ for some $\eta\ge0$ then there exists $\theta$ depending on $p, q, \eta, n$, $w$ and $\kappa$ {\upshape(}the constant from \eqref{slow growth}{\upshape)} and a constant $C>0$ depending on $\theta$ such that 
\begin{align*}
	\Big\Vert \Big(\sum_j \big|\MM^\theta(f_j)\big|^q\Big)^{1/q}\Big\Vert_{L^p_w} \le C\Big\Vert\Big(\sum_j |f_j|^q\Big)^{1/q}\Big\Vert_{L^p_w}.
\end{align*}

\item If $w\in \AH_\infty$ then for any $0<s<\min\{p/r_w,q\}$, there exists $\theta$ depending on $p, q, s, n$, $w$ and $\kappa$ {\upshape(}the constant from \eqref{slow growth}{\upshape)} and a constant $C>0$ depending on $\theta$ such that
\begin{align*}
	\Big\Vert \Big(\sum_j \big|\MM^\theta_s(f_j)\big|^q\Big)^{1/q}\Big\Vert_{L^p_w} \le C\Big\Vert\Big(\sum_j |f_j|^q\Big)^{1/q}\Big\Vert_{L^p_w}.
\end{align*}

\end{enumerate} 
\end{Lemma}
\begin{proof}
The proof of part (a) can be found in the proof of \cite[Theorem 2.7]{Tang}. To prove (b) we first observe that since $w\in\AH_\infty$ then $w\in\AH_\infty^\eta$ for some $\eta\ge 0$. Now if $s<\f{p}{r_w}$ then $w\in\AH_{p/s}^\eta$ by Remark \ref{rem: weights}. We may then apply part (a) because $p/s>1$ and $q/s>1$. Thus we have
	\begin{align*}
		\Big\Vert \Big(\sum_j \big|\MM_s^\theta(f_j)\big|^q\Big)^{\f{1}{q}}\Big\Vert_{L^p_w} 
		\lesssim   \Big\Vert \Big(\sum_j \big| f_j^s \big|^{\f{q}{s}}\Big)^{\f{s}{q}}\Big\Vert^{\f{1}{s}}_{L^{p/s}_w}
		=   \Big\Vert \Big(\sum_j \big| f_j \big|^{q}\Big)^{\f{1}{q}}\Big\Vert_{L^{p}_w}. \qquad\qedhere
	\end{align*}
\end{proof}

\subsection{Hermite tiles and sequence spaces}\label{sec:tiles}
The aim of this section is to detail an appropriate notion of discretized sets for $\RR^n$ and then use these sets to define weighted sequence spaces in Definition \ref{def:sequences}.

It is well known that the dyadic cubes of $\RR^n$ play an important role in the development of frames for the classical distribution spaces. In the Hermite setting, the notion of `tiles' or `rectangles' constructed in \cite{PX} will play the role of dyadic cubes for $\LL$. The construction of these tiles rely on the zeros of the Hermite polynomials. 

Fix a constant $\delta_\star\in(0,\f{1}{37})$ and for each $j\ge 0$ we set
$$ N_j = \big[(1+11\delta_\star)(\tfrac{4}{\pi})^2 4^j\big] +3. $$
Let $\{\zeta_\nu\}$ with $\nu\in \{\pm 1, \pm 2, \dots, \pm N_j\}$ be the zeros of the $H_{2N_j}(t)$, ordered so that
\begin{align}\label{zeros}
 \zeta_{-N_j}<\dots<\zeta_{-1}<0<\zeta_1<\dots<\zeta_{N_j}.
 \end{align}
It will be useful to note that $\zeta_{-\nu} =- \zeta_\nu$.

The index sets $\X_j$ are formed from the zeroes \eqref{zeros} of the even Hermite polynomials as follows. For each $j\ge 0$ we set 
$$ \X_j := \big\{ (\zeta_{\alpha_1},\dots, \zeta_{\alpha_n}) \in \RR^n: |\alpha_\nu|<N_j, \;\nu=1,\dots,n \big\}.$$
Following the convention in \cite{PX} the elements of $\X_j$ will be called \emph{nodes}. We write $\X=~\bigcup_{j\ge 0}\X_j$ to represent the collection of all nodes.

The families of \emph{dyadic tiles} $\E_j$ for $j\ge 0$ are defined as follows. For each 
node $\zeta\in \X_j$  with $\zeta = (\zeta_{\alpha_1},\dots, \zeta_{\alpha_n})$ we set 
$$ R_{\zeta}:= I_{\alpha_1}\times I_{\alpha_2}\times \dots \times I_{\alpha_n},  $$
where $I_{\alpha_\nu}$ are the intervals defined by
	\begin{align*}
		I_{\nu}:= \left
\lbrace 
	\begin{array}{ll}
	\left[0, (\zeta_1+\zeta_2)/2\right], &\qquad \nu=1\\
	\left[(\zeta_{\nu-1}+\zeta_\nu)/2, (\zeta_\nu+\zeta_{\nu+1})/2\right], &\qquad \nu \in \{ 2,3,\dots, N_{j-1}\}\\
	\left[(\zeta_{N_{j-1}}+\zeta_{N_j})/2, \zeta_{N_j}+2^{-j/6}\right], &\qquad \nu=N_j
	\end{array}
\right.
\end{align*}
and $I_{-\nu}=-I_{\nu}$ for $\nu\in\{1,2\dots,N_j\}$.

We denote by $\E_j$ the collection of tiles $\{R_\zeta\}_{\zeta\in \X_j}$ and by $\E := \bigcup_{j\ge 0}\E_j$ the collection of all tiles. Then $\E_j$ contains approximately $4^{jn}$ tiles (with constants depending on $\delta_\star$ and $n$).  
By a rectangle (or tile) $R$ we mean a tile from $\E_j$ for some $j\ge 0$ and denote its node by $x_R$.
We summarize the properties of these tiles below. For their proofs we direct the reader to \cite[pp. 379-380]{PX}.

\begin{Lemma}\label{lem: tiles}
There exist constants $c_0$, $c_1$, $c_2$, $c_3$ and $c_4$ depending only on $\delta_\star$ and $n$ such that for each $j\ge 0$ and each tile $R\in\E_j$ the following holds. 
\begin{enumerate}[\upshape(a)]
\item If $|x_R|\le (1+4\delta_\star)2^{j+1}$, it holds that
\begin{align}\label{Hermite tiles 0} 
R\subset Q(x_R,c_0 2^{-j}). 
\end{align}
\item In general, it holds that
\begin{align}\label{Hermite tiles 1-2}
 Q(x_R, c_1 2^{-j})\subset R \subset Q(x_R, c_2 2^{-j/3}).
 \end{align}
\item Set $\Q_j:=\bigcup\limits_{P\in\E_j}P = Q(0,\zeta_{N_j}+2^{-j/6})$; it holds that
\begin{align}\label{Hermite tiles 3} Q(0,2^j)\subset \Q_j \subset Q(0,c_32^j).
\end{align}
\item $R$ can be subdivided into a disjoint union of subcubes of sidelength roughly equal to $2^{-j}$; more precisely each such subcube $Q$ satisfies
$$ Q(x_Q, c_4 2^{-j-1})\subset Q \subset Q(x_Q, c_4 2^{-j}). $$
Denoting by $\widehat{\E}_j$ the collection of all subdivided cubes from $\E_j$, then it holds that 
\begin{align}\label{Hermite cubes}
 \bigcup_{Q\in\widehat{\E}_j}Q= \Q_j.
 \end{align}
 \end{enumerate}
 \end{Lemma}

We may now introduce our weighted sequence spaces associated with Hermite expansions, which are spaces of sequences indexed over the collection of tiles just described. 

\begin{Definition}[Weighted Hermite sequence spaces]\label{def:sequences}
Let  $\alpha\in\RR$, $0<q\le\infty$  and $0\le w\in L^1_{loc}(\RR^n)$.
For $0<p\le\infty,$  we define the weighted Hermite Besov sequence space $b^{p,q}_{\alpha,w}$ as the set of all sequences of complex numbers $s=\{s_R\}_{R\in \E}$ such that 
$$ \Vert s\Vert_{b^{p,q}_{\alpha, w}} :=\bigg\{ \sum_{j\ge 0} 2^{j\alpha q}\Big(\sum_{R\in\E_j} \big(w(R)^{1/p} |R|^{-1/2}|s_R|\big)^p \Big)^{q/p}\bigg\}^{1/q} <\infty;
$$
for $0<p<\infty,$  we define the  weighted Hermite Triebel--Lizorkin sequence space $f^{p,q}_{\alpha,w}$ as the set of all sequences of complex numbers $s=\{s_R\}_{R\in \E}$ such that
$$ \Vert s\Vert_{f^{p,q}_{\alpha,w}}:= \bigg\Vert \Big(\sum_{j\ge 0} 2^{j\alpha q}\sum_{R\in\E_j}\big(\Ind_R(\cdot)|R|^{-1/2} |s_R|\big)^q\Big)^{1/q}\bigg\Vert_{L^p_w} <\infty.
$$
We use the supremum norm $\ell^\infty$ when $q=\infty$ (or if $p=\infty$ in the Besov case).
\end{Definition}
These spaces are critical for our frame characterizations in Section \ref{sec:frames} and the embedding results of Section \ref{sec:embedding}.

Throughout the rest of this article, we will use the notation $A^{p,q}_{\alpha,w}(\LL)$  (or $A^{p,q}_{\alpha,w}$) to refer to $B^{p,q}_{\alpha,w}(\LL)$ or $F^{p,q}_{\alpha,w}(\LL),$ with the understanding that  $\alpha\in \RR,$ $0<q\le \infty,$ and either $0<p\le \infty$ if $A=B$ or $0<p<\infty$ if $A=F.$ An analogous comment applies to the sequence spaces, denoted by $a^{p,q}_{\alpha,w}(\LL)$ (or $a^{p,q}_{\alpha,w}$).

\section{Frame characterizations}\label{sec:frames}
In this section we obtain frame characterizations of our weighted Hermite spaces. 
Frames  were constructed in \cite{PX} for the unweighted setting and here we show that the same notions can be used to characterize the weighted spaces. Our frames rely on a certain cubature formula for functions in $V_N$ (the spaces of projections defined in Section \ref{sec:prelim}). Consider the well known  `Christoffel function' 
$$ \tau(N,x):= \f{1}{\QQ_N(x,x)},\qquad x\in \RR, \;N\in\NN_0, $$
which has certain useful asymptotic properties (see \cite[p. 376]{PX}). Then  the following cubature formula 
	\begin{align}\label{eq:cubature}
	 \int_{\RR^n} f(x)\,g(x)\,dx \sim \sum_{\zeta\in \X_j} \tau_\zeta\, f(\zeta)\,g(\zeta),\qquad \zeta = (\zeta_{\alpha_1},\dots,\zeta_{\alpha_n}), \;\tau_\zeta = \prod_{k=1}^n \tau(2N_j,\zeta_{\alpha_k}) \end{align}
	is exact for all $f\in V_k$ and $g\in V_\ell$ with $k+\ell \le 4N_j-1$. See \cite[Proposition 2]{PX}.

If $\{\varphi_j\}_{j\ge 0}$ is an admissible system then for  each tile $R\in \E_j$ we set
\begin{align}\label{eq:needlet} \varphi_R(x):=\tau_R^{1/2} \varphi_j(\sqrt{\LL})(x,x_R)\end{align}
where $x_R$ is the node of $R$ and $\tau_R = \tau_{x_R}$, the coefficient in the cubature formula \eqref{eq:cubature}. They satisfy $|\tau_{R}|\sim |R|$ for any tile $R$ (see \cite[(2.33)]{PX}). The functions $\vph_R$ were termed \emph{needlets} in \cite{PX}. 

Given any admissible systems $\{\varphi_j\}_{j\ge 0}$ and $\{\psi_j\}_{j\ge 0}$ we define the  \emph{analysis} $S_\varphi$ and \emph{synthesis} $T_\psi$ operators by
\begin{align*}
	S_\varphi: f\longmapsto \{ \ip{f,\varphi_R}\}_{R\in\E} &&\text{and}&&
	T_\psi: \{s_R\}_{R\in\E}\longmapsto \sum_{R\in\E} s_R\psi_R.
\end{align*}
The main result in this section is the following. For the unweighted case see \cite[Theorems~3~\&~5]{PX}.
\begin{Theorem}[Frame characterization]\label{thm:frames}
Let $w\in\AH_\infty$,  $\alpha\in\RR$, $0<q\le \infty,$ and $0<p<\infty$ if $A^{p,q}_{\alpha,w}(\LL)=F^{p,q}_{\alpha,w}(\LL)$ or $0<p\le\infty$ if $A^{p,q}_{\alpha,w}(\LL)=B^{p,q}_{\alpha,w}(\LL)$. Suppose that $\{\varphi_j\}_{j\ge 0}$ and $\{\psi_j\}_{j\ge 0}$ are two admissible systems. Then,
\begin{enumerate}[\upshape(a)]
	\item the operator $T_\psi: a^{p,q}_{\alpha,w}(\LL)\to A^{p,q}_{\alpha,w}(\LL)$ is bounded;
	\item the operator $S_\varphi: A^{p,q}_{\alpha,w}(\LL)\to a^{p,q}_{\alpha,w}(\LL)$ is bounded;
	\item if $\{\varphi_j\}_{j\ge 0}$ and $\{\psi_j\}_{j\ge 0}$   satisfy
	\begin{align}\label{eq:CRF1}
	\sum_{j\ge 0}\psi_j(\lambda)\varphi_j(\lambda) = 1 \qquad \forall \lambda\ge 0,
	\end{align}
	 then $T_\psi\circ S_\varphi = I$ on $A^{p,q}_{\alpha,w}(\LL)$, with convergence in $\sz'(\RR^n)$. 
	 Furthermore, the definitions of $A^{p,q}_{\alpha,w}(\LL)$ are independent of the choice of $\{\varphi_j\}_{j\ge 0}$.
\end{enumerate}
\end{Theorem}
The proof of Theorem \ref{thm:frames} is given in Section \ref{sec:frameproof}. Before embarking on the proof we need to give some preparatory facts and estimates. Our first set of results deals with the analysis and synthesis operators, which adapts \cite[Proposition 3]{PX} to our band-limited operators.

\begin{Proposition}\label{prop: CRF}
Let $\{\varphi_j\}_{j\ge 0}$ and $\{\psi_j\}_{j\ge 0}$ be two admissible systems satisfying \eqref{eq:CRF1}.
Then we have the following:
	\begin{align}\label{eq:CRF2} f = \sum_{j\ge 0} \psi_j(\sqrt{\LL})\varphi_j(\sqrt{\LL}) f \qquad \text{in}\qquad  \sz'(\RR^n) \end{align}
and
	\begin{align}\label{eq:CRF3} f = \sum_{R\in \E} \ip{f, \varphi_R} \,\psi_R \qquad \text{in}\qquad  \sz'(\RR^n). \end{align}
\end{Proposition}
\begin{Remark}\label{rem:CRF}
It is worth pointing out that if \eqref{eq:CRF1} holds then Proposition \ref{prop: CRF} implies that  $T_\psi\circ S_\varphi = I$ on $\sz'(\RR^n)$. Another consequence of Proposition \ref{prop: CRF} is that our spaces $B^{p,q}_{\alpha,w}(\LL)$ and $F^{p,q}_{\alpha,w}(\LL)$ are quasi-Banach spaces embedded continuously into $\mathscr{S}'(\RR^n)$. To see this one can reason as in \cite[Proposition 4]{PX} and \cite[Section 5.1]{PX}.
\end{Remark}
\begin{proof}[Proof of Proposition \ref{prop: CRF}]
The reproducing formula \eqref{eq:CRF2} can be found in \cite[Proposition~5.5(b)]{KP} and \cite[Proposition 3]{PX}, so we only address \eqref{eq:CRF3}.

To prove \eqref{eq:CRF3} we firstly have $\varphi_j(\sqrt{\LL})(x,\cdot)\in V_{ 4^j}$. Indeed from \eqref{eq:kernel} we see that $\vph_j(\sqrt{\LL})(x,\cdot)$ is a finite combination of Hermite functions $h_\mu$, of order $|\mu|$ not exceeding $4^j$ which, by the definition of $V_N$, implies $\varphi_j(\sqrt{\LL})(x,\cdot)\in V_{4^j}$. One can reason similarly for $\psi_j(\sqrt{\LL})(x,\cdot)$. 
Therefore since $ 2\cdot4^{j}\le 4 N_j -1$ we may apply \eqref{eq:cubature} to $\varphi_j(\sqrt{\LL})(x,\cdot)$ and $\psi_j(\sqrt{\LL})(x,\cdot)$ to obtain
\begin{align*}
	\psi_j(\sqrt{\LL})\varphi_j(\sqrt{\LL})(x,y) 
	= \int \psi_j(\sqrt{\LL})(x,u)\,\varphi_j(\sqrt{\LL})(u,y)\,du
	= \sum_{\zeta \in\X_j} \tau_\zeta \,\psi_j(\sqrt{\LL})(x,\zeta)\,\varphi_j(\sqrt{\LL})(\zeta,y).
\end{align*}
Then in view of \eqref{eq:needlet} we see that
\begin{align*}
	\psi_j(\sqrt{\LL})\varphi_j(\sqrt{\LL})(x,y) 
	=\sum_{R\in\E_j} \tau_R^{1/2} \psi_j(\sqrt{\LL})(x,\xi)\,\tau_R^{1/2} \varphi_j(\sqrt{\LL})(y,\xi)
	=\sum_{R\in\E_j} \psi_R(x)\,\varphi_R(y).
\end{align*}
The result now follows by invoking \eqref{eq:CRF2} and then utilising the above equality. Thus we have
\begin{align*}
	f(x) 
	=\sum_{j\ge 0} \int \psi_j(\sqrt{\LL})\varphi_j(\sqrt{\LL}) (x,y)f(y)\,dy
	=\sum_{j\ge 0}\sum_{R\in\E_j} \Big(\int \varphi_R(y) f(y)\,dy\Big) \psi_R(x),
\end{align*}
which is \eqref{eq:CRF3}.
\end{proof}

The next estimate will be needed in the proofs of Theorems \ref{thm:frames} and \ref{thm:embedding2}.
\begin{Lemma}\label{almost orthog}
Let $\{\varphi_j\}_{j\ge 0}$, $\{\psi_j\}_{j\ge 0}$ be any pair of admissible systems and let $\sigma\ge1$. Then there exists $C>0$ depending on $\sigma$ such that for any sequence of numbers $s=\{s_R\}_{R\in\E}$,
\begin{align*}
		\big| \varphi_j(\sqrt{\LL})(T_\psi s)(x)\big| \le C \sum_{k=j-2}^{j+2}\sum_{R\in\E_k} \f{|R|^{-1/2}|s_R|}{(1+2^k|x-x_R|)^\sigma}, \qquad \forall x\in\RR^n.
\end{align*}
\end{Lemma}
\begin{proof}[Proof of Lemma \ref{almost orthog}]
Recall from \eqref{eq:orthog1} that $\varphi_j(\sqrt{\LL})\psi_R(x)=0$, whenever $R\in \E_k$ with $|j-~k|\ge 3$. 
It follows that
	\begin{align*}
		\varphi_j(\sqrt{\LL}) (T_\psi s)(x)
		=\sum_{k=j-2}^{j+2} \sum_{R\in\E_k} s_R \,\tau_R^{1/2} \varphi_j(\sqrt{\LL})\psi_k(\sqrt{\LL})(x,x_R).
	\end{align*}
	Then by \eqref{eq:orthog2} with $\ve =4(1+4\delta_\star)^2$ we obtain for any $\sigma\ge 1$,
	\begin{align*}
		\big|\varphi_j(\sqrt{\LL}) (T_\psi s)(x)\big|
\lesssim \sum_{k=j-2}^{j+2} \sum_{R\in\E_k}  |s_R| \,|R|^{1/2}\f{2^{kn}\e_{\ve 4^k}(x_R)\e_{\ve 4^{j}}(x) }{(1+2^k|x-x_R|)^\sigma}.
	\end{align*}
The proof is concluded by invoking the following estimate:
\begin{align*}
2^{kn} \e_{\ve 4^k}(x_R)\lesssim |R|^{-1}\quad \forall\; R\in \E_k, k\in \NN_0.
\end{align*}	
We can see this estimate by considering the two types of tiles listed in Lemma~\ref{lem: tiles}. Indeed if $|x_R|_\infty \le (1+4\delta_\star)2^{k+1}=\sqrt{\ve 4^k}$, then $2^{kn}\sim |R|^{-1}$ by \eqref{Hermite tiles 0}. On the other hand, if $|x_R|_\infty > (1+4\delta_\star)2^{k+1}$, then by  \eqref{eq-eN} we have $\e_{\ve 4^k}(x_R)=e^{-\vartheta |x_R|^2} \le e^{-\vartheta 4^k} \lesssim 2^{-\beta k}$ for any $\beta>0$. Taking $\beta=2n/3$ in view of the inequality $2^{kn/3}\lesssim |R|^{-1}$ (see \eqref{Hermite tiles 1-2}) gives the required result.
\end{proof}
\subsection{Maximal lemmas} 
In this section we gather and prove some maximal lemmas  for sequences (Lemmas \ref{lem:peetre}, \ref{lem: seq} and \ref{lem:plancherel-polya}) that will be needed in the proof of Theorem~\ref{thm:frames}. Recall that the maximal function $\MM_s^\theta$ was defined in \eqref{mf}.

We first have a generalization of \cite[Lemma 4]{PX}.
\begin{Lemma}\label{lem:peetre}
Let $s>0$, $\theta\ge 0$, $\sigma>\f{\theta}{s}+\max\{n,\f{n}{s}\}$ and $j\ge 0$.
For a collection of numbers $\{a_R\}_{R\in\E_j}$ we set
\begin{align}\label{eq:sequence sup}
	a^*_j(x):=\sum_{R\in\E_j}\f{|a_R|}{(1+2^j|x-x_R|)^\sigma}, \qquad x\in\RR^n,
\end{align}
and $a^*_R:=a^*_j(x_R)$.
Then there exists $C>0$ such that for any $\tilde{c}\in(0,2c_4]$ and $x\in\RR^n$ we have
\begin{align}\label{eq:maximal A}
a^*_j(x)\le C \inf_{y\in Q(x,\tilde{c}2^{-j})} \MM^\theta_s\Big(\sum_{R\in\E_j}|a_R|\Ind_R\Big)(y),
\end{align}
and
\begin{align}\label{eq:maximal B}
a^*_R\Ind_R(x)\le C \inf_{y\in Q(x,\tilde{c}2^{-j})} \MM^\theta_s\Big(\sum_{R\in\E_j}|a_R|\Ind_R\Big)(y).
\end{align}
Here the constants depend only on $n$, $\delta_\star$, $\sigma$, $\theta$ and $s$.
\end{Lemma}
\begin{proof}[Proof of Lemma \ref{lem:peetre}]

We first set 
\begin{align*}
	\widetilde{a}_j(x):=\sum_{R\in \E_j} \f{|a_R|}{(1+2^jd(x,R))^\sigma}, 
\end{align*}
where $d(x,R):=\sup\limits_{y\in R}|x-y|_\infty$. Since $a^*_j(x)\lesssim \tilde{a}_j(x)$ and $a^*_R\Ind_R(x)\lesssim \tilde{a}_j(x)$ 
for every $x\in\RR^n$, in order to obtain \eqref{eq:maximal A} and \eqref{eq:maximal B} it suffices to prove 
\begin{align}\label{maximal 1}
\widetilde{a}_j(x)\lesssim \inf_{y\in Q(x,\tilde{c}2^{-j})} \MM^\theta_s\Big(\sum_{R\in\E_j}|a_R|\Ind_R\Big)(y), \qquad \forall x\in\RR^n.
\end{align}
Set $\nu=1-\min\{1,1/s\}$ and $\widetilde{Q}:=Q(x,\tilde{c}2^{-j})$. Let $c_3$ and $c_4$ be the constants from Lemma~\ref{lem: tiles}. We split the proof of \eqref{maximal 1} into two cases. 

\underline{Case 1: $|x|_\infty>2(c_3+c_4)2^j$.}

For each such $x$ we have $d(x,R)>\f{|x|_\infty}{2}$ for any $R\in\E_j$. Then by either H\"older's inequality if $s>1$ or the $s$-triangle inequality if $s\le 1$, we have
\begin{align*}
	\widetilde{a}_j(x)
	\lesssim \Big(\f{2^{-j}}{|x|_\infty}\Big)^\sigma \sum_{R\in\E_j}|a_R|
	\lesssim \f{2^{-j(\sigma-2n\nu)}}{|x|_\infty^\sigma}  \Big(\sum_{R\in\E_j}|a_R|^s\Big)^{1/s}.
\end{align*}
Set $Q_x:=Q(0,2|x|_\infty)$. Then Lemma \ref{lem: tiles} (c) ensures $\Q_j\subseteq Q_x$. 
Now using the $1/s$-triangle inequality if $s\ge 1$, and  H\"older's inequality otherwise we have
\begin{align*}
\Big(\sum_{R\in\E_j}|a_R|^s\Big)^{1/s}
=\Big(\int_{\Q_j}\sum_{R\in\E_j} |a_R|^s\f{\Ind_R(y)}{|R|}\,dy\Big)^{1/s}
\le \Big(\int_{\Q_j}\Big(\sum_{R\in\E_j} |a_R|\f{\Ind_R(y)}{|R|^{1/s}}\Big)^s\,dy\Big)^{1/s}.
\end{align*}
Next using the estimates $|Q_x|\sim |x|_\infty^n$, $\Psi_\theta(Q_x)\sim |x|^\theta_\infty$ and  $|R| \gtrsim 2^{-jn}$, valid for every $R\in\E_j$, we obtain
\begin{align*}
\int_{\Q_j}\Big(\sum_{R\in\E_j} |a_R||R|^{-1/s}\Ind_R(y)\Big)^s\,dy
	&\lesssim  2^{jn} |x|_\infty^{n+\theta}\f{1}{\Psi_\theta(Q_x)}\aver_{Q_x}\Big(\sum_{R\in\E_j} |a_R| \Ind_R(y)\Big)^s\,dy.
\end{align*}
Combining the previous three estimates and observing that $\sigma > \f{n+\theta}{s}$ we have
\begin{align*}
\widetilde{a}_j(x) \lesi 2^{-j(2\sigma-\f{2n}{s}-\f{\theta}{s}-2n\nu)}\inf_{y\in \widetilde{Q}}\MM_s^\theta\Big(\sum_{R\in\E_j}|a_R|\Ind_R\Big)(y).
\end{align*}
We arrive at \eqref{maximal 1} on recognizing that $\sigma > \f{n}{s}+\f{\theta}{2s}$ if $s\le 1$ and $\sigma>n+\f{\theta}{2s}$ if $s\ge 1$, both of which hold because of our assumption on $\sigma$.

\medskip
\underline{Case 2: $|x|_\infty\le2(c_3+c_4)2^j$.}

Let $\widehat{\E}_j$ be the collection of cubes defined in \eqref{Hermite cubes}. Then for each $Q\in\widehat{\E}_j$, we set $a_Q:=a_R$ whenever $Q\subset R$. Then we have
\begin{align}\label{maximal 2}
	\widetilde{a}_j(x)\le \sum_{Q\in\widehat{\E}_j} \f{|a_Q|}{(1+2^jd(x,Q))^\sigma}
	\quad \text{and}\quad
	\sum_{R\in\E_j}|a_R|\Ind_R = \sum_{Q\in\widehat{\E}_j}|a_Q|\Ind_Q.
\end{align}
We shall discretize $\RR^n$ into the following `square annuli' and cube centred at $x$. For each $m\ge 1$, we set
\begin{align*}
\A_0&=\A_0(x,j)=\{Q\in\widehat{E}_j: |x-x_Q|_\infty\le c_4 2^{-j}\}, \\
\text{and} \ \ \A_m &= \A_m(x,j)=\{Q\in\widehat{E}_j: c_4 2^{m-j-1}<|x-x_Q|_\infty\le  c_4 2^{m-j}\}.
\end{align*}
We shall also need the the following cubes centred at $x$. For each $m\ge 0$, we set
\begin{align*}
\B_m&=\B_m(x,j)
= Q(x,c_42^{m-j+1}).
\end{align*}
These sets satisfy the following properties
\begin{align}\label{eq:maximal3}
 \# \A_m \lesssim 2^{mn}, 
 &&
 \widehat{\E}_j = \bigcup\limits_{m\ge 0} \A_m,
 &&
 \bigcup_{Q\in\A_m} Q\subset \B_m.
 \end{align}
By \eqref{maximal 2} and H\"older's inequality if $s>1$ and the $s$-triangle inequality otherwise,
 \begin{align*}
 	\widetilde{a}_j(x)
	\le \sum_{m\ge 0}\sum_{Q\in\A_m} \f{|a_Q|}{(1+2^jd(x,Q))^\sigma} 
	\lesssim \sum_{m\ge 0} 2^{-m(\sigma-n\nu)} \Big(\sum_{Q\in\A_m} |a_Q|^s\Big)^{1/s}.
\end{align*}
Using the last property in  \eqref{eq:maximal3}, along with the $1/s$-triangle inequality if $s\ge 1$ and H\"older's inequality otherwise, we have
 \begin{align*}
\Big(\sum_{Q\in\A_m} |a_Q|^s\Big)^{1/s}
\le  \Big(\int_{\B_m} \sum_{Q\in\A_m} |a_Q|^s |Q|^{-1}\Ind_Q(y)\,dy \Big)^{1/s}.
\end{align*}
Now using the estimates $|\B_m|\sim 2^{n(m-j)}$, $|Q|\sim 2^{-jn}$, $\Psi_\theta(\B_m)\lesssim 2^{(m-j)\theta}|x|_\infty^\theta$ and  $|x|_\infty\lesssim 2^j$ with the fact that $\B_m$ contains $x$ we have
  \begin{align*}
  \int_{\B_m} \sum_{Q\in\A_m} |a_Q|^s |Q|^{-1}\Ind_Q(y)\,dy
\lesi 2^{m(n+\theta)} \f{1}{\Psi_\theta(\B_m)}\aver_{\B_m} \Big(\sum_{Q\in\A_m} |a_Q| \Ind_Q(y)\Big)^s\,dy.
 \end{align*}
Combining the previous three calculations gives
  \begin{align*}
\widetilde{a}_j(x) \lesi  \inf_{y\in \widetilde{Q}}\MM_s^\theta\Big(\sum_{R\in\E_j}|a_R|\Ind_R\Big)(y) \times\sum_{m\ge 0} 2^{-m(\sigma-n\nu-\f{n}{s}-\f{\theta}{s})}.
 \end{align*}
 Our assumption on $\sigma$ ensures the convergence of the sum,
giving \eqref{maximal 1} and completing the proof of Lemma~\ref{lem:peetre}.
\end{proof}

\begin{Lemma}[{\cite[Lemma 5]{PX}}]\label{lem: seq}
Fix $j\ge0$ and let $g\in V_{4^j}$. For each tile $R\in\E_j$ define
\begin{align}\label{seq1}
	g_R=\sup_{x\in R}|g(x)| && \text{and} &&  \mathfrak{g}_R=\inf_{x\in R}|g(x)|.
\end{align}
Then there exist $k\ge 1$ and $C>0$ depending only on $n, \delta_\star$ and $\sigma$ such that for every $R\in \E_j$ we have
\begin{align*}
	g_R^* \le C  \mathfrak{g}_P^*, \qquad \forall P\in \E_{j+k},\quad\text{with}\quad P\cap R\ne \emptyset,
\end{align*}
and
\begin{align}\label{eq:seq3}
	g_R^*\Ind_R(x)\le C \sum_{\substack{P\in\E_{j+k}, \\ P\cap R\ne\emptyset}} \mathfrak{g}_P^*\Ind_P(x),\qquad\forall x\in \RR^n.
\end{align}
\end{Lemma}

\begin{Lemma}\label{lem:plancherel-polya}
Fix $0<p\le \infty$ and $w\in\AH_\infty$. Then there exists $C>0$ such that for any $j\ge 0$ and $g\in V_{4^j}$ we have
$$ \Big(\sum_{R\in\E_j} w(R) \max_{x\in R} | g(x)|^p\Big)^{1/p}\le C \Vert g\Vert_{L^p_w}.$$
\end{Lemma}
\begin{proof}[Proof of Lemma \ref{lem:plancherel-polya}]
For any tile $R\in \E_j$ let $g_R$ and $\mathfrak{g}_R$ denote the collection of numbers from \eqref{seq1}. Let $k$ be the integer from Lemma \ref{lem: seq}.   Fix  $s\in (0, \min\{p/r_w, 1\})$ and  let $\theta$ be the number from Lemma \ref{lem:FS}.

From the fact that the tiles in $\E_j$ are all disjoint, then the inequality $a_R\le a_R^*$ and \eqref{eq:seq3} give
\begin{align*}
	\Big(\sum_{R\in\E_j} w(R)\max_{x\in R} | g(x)|^p\Big)^{1/p}
	\le \bigg\Vert\sum_{R\in\E_j} g_R^*\Ind_R\bigg\Vert_{L^p_w}
	\lesssim \bigg\Vert \sum_{R\in\E_j}\sum_{\substack{P\in\E_{j+k}, \\ P\cap R\ne\emptyset}} \mathfrak{g}_P^*\Ind_P\bigg\Vert_{L^p_w}.
\end{align*}
Now note that $\bigcup_{R\in\E_j}R\subset \bigcup_{P\in\E_{j+k}}P$. Applying \eqref{eq:maximal B} with $\sigma>\theta/s + \max\{n,\f{n}{s}\}$ we have
\begin{align*}
\bigg\Vert \sum_{R\in\E_j}\sum_{\substack{P\in\E_{j+k}, \\ P\cap R\ne\emptyset}} \mathfrak{g}_P^*\Ind_P\bigg\Vert_{L^p_w}
\le  \bigg\Vert \sum_{P\in\E_{j+k}} \mathfrak{g}_P^*\Ind_P\bigg\Vert_{L^p_w} 
	\lesssim \bigg\Vert \MM_s^\theta\Big(\sum_{P\in\E_{j+k}} \mathfrak{g}_P\Ind_P\Big)\bigg\Vert_{L^p_w}.
\end{align*}
Finally we invoke Lemma \ref{lem:FS} (b) with $q=1$ to obtain
\begin{align*}
\Big(\sum_{R\in\E_j} w(R)\max_{x\in R} | g(x)|^p\Big)^{1/p}
\lesi \bigg\Vert \sum_{P\in\E_{j+k}} \mathfrak{g}_P\Ind_P\bigg\Vert_{L^p_w} 
	\le \Vert g\Vert_{L^p_w},
\end{align*}
which concludes the proof.
\end{proof}

\subsection{Proof of Theorem \ref{thm:frames}}\label{sec:frameproof}
In this section we  prove the frame characterizations for our weighted Besov and Triebel--Lizorkin spaces.  We start by observing that $T_\psi\circ S_\vph=I$ follows readily from  Proposition \ref{prop: CRF} (see Remark \ref{rem:CRF}). Furthermore, the independence of $\vph$ in the definitions of the spaces $A^{p,q}_{\alpha,w}(\LL)$ can be seen  by following a similar argument to the unweighted case (see \cite[Theorems 3 \& 5]{PX}). This requires proving part (a) assuming that the spaces $A^{p,q}_{\alpha,w}$ have been defined using some other admissible pair $(\wt{\vph}_\circ,\wt{\vph})$  different from the one used in the definition of $S_\vph$. 

Thus we only need to prove (a) and (b), and we shall address the Triebel--Lizorkin and Besov scales  separately.

\medskip
\underline{\it The Triebel--Lizorkin case.}
Part (a). Suppose that $F^{p,q}_{\alpha,w}$ has been defined using some admissible pair $(\wt{\vph}_\circ,\wt{\vph})$. 
Let $0<r<\min\{ p/r_w,q\}$ and $\theta$ be the number in Lemma \ref{lem:FS} (b). Fix $\sigma> \theta/r +\max\{n, n/r\}$. By Lemma \ref{almost orthog} and \eqref{eq:maximal A} with $a_R= |R|^{-1/2}|s_R|$ we have
\begin{align*}
\big|\wt{\varphi}_j(\sqrt{\LL})(T_{\psi} s)(x)\big|
	\lesssim \sum_{k=j-2}^{j+2}\sum_{R\in\E_k} \f{|R|^{-1/2}|s_R|}{(1+2^k|x-x_R|)^\sigma}
	\lesssim \sum_{k=j-2}^{j+2}\MM^\theta_r\Big(\sum_{R\in\E_k} |R|^{-1/2}|s_R| \Ind_R\Big)(x).
\end{align*}
Applying this estimate along with Lemma \ref{lem:FS} (b) we have
\begin{align*}
	\Vert T_{\psi} s\Vert_{F^{p,q}_{\alpha,w}}
	&\lesssim \Big\Vert \Big(\sum_{j\ge 0}\Big(2^{j\alpha}  \sum_{k=j-2}^{j+2}\MM^\theta_r\Big(\sum_{R\in\E_k} |R|^{-1/2}|s_R| \Ind_R\Big)\Big)^q \Big)^{1/q}\Big\Vert_{L^p_w}\\
	&\lesssim \Big\Vert \Big(\sum_{j\ge 0} \MM^\theta_r\Big(2^{j\alpha}  \sum_{R\in\E_j} |R|^{-1/2}|s_R| \Ind_R\Big)^q \Big)^{1/q}\Big\Vert_{L^p_w}\\
	&\lesssim \Big\Vert \Big(\sum_{j\ge 0} \Big(2^{j\alpha}  \sum_{R\in\E_j} |R|^{-1/2}|s_R| \Ind_R\Big)^q \Big)^{1/q}\Big\Vert_{L^p_w}
	= \Vert s\Vert_{f^{p,q}_{\alpha,w}}.
\end{align*}

Part (b). Suppose now that $F^{p,q}_{\alpha,w}$ has been defined using the admissible pair $(\vph_\circ,\vph)$ from the definition of $S_\vph$. Let $f\in F^{p,q}_{\alpha, w}$. Note that for each $j\ge 0$ we have $\varphi_j(\sqrt{\LL})f\in V_{ 4^j}$ and hence Lemma \ref{lem: seq} can be applied to $g=\varphi_j(\sqrt{\LL})f$. 
Define the sequences 
\begin{align*}
	\ip{f^{\varphi, R}} &:= \sup_{x\in R} | \varphi_j(\sqrt{\LL}) f(x)|, \qquad R\in \E_j, \\
	\ip{f_{\varphi, P, k}}&:=\inf_{x\in P} |\varphi_j(\sqrt{\LL})f(x)|,\qquad P\in \E_{j+k}.
\end{align*}	
Then we firstly have
\begin{align}\label{ff1}
	|\ip{f,\varphi_R}| = \big| \tau_R^{1/2} \varphi_j(\sqrt{\LL})f(x_R)\big| 
	\lesssim |R|^{1/2} \ip{f^{\varphi,R} }
	\le |R|^{1/2} \ip{f^{\varphi,R}}^*,
\end{align}
where $\ip{f^{\varphi,R}}^*$ is defined in \eqref{eq:sequence sup}. Secondly, we apply Lemma \ref{lem: seq} to $g_R=\ip{f^{\varphi,R}}$ and $\mathfrak{g}_P = \ip{f_{\varphi,P,k}}$ to obtain
\begin{align}\label{ff2}
	\ip{f^{\varphi,R}}^* \Ind_R(x)\lesssim \sum_{\substack{P\in\E_{j+k}, \\ P\cap R\ne\emptyset}} \ip{f_{\varphi,P,k}}^*\Ind_P(x), \qquad x\in\RR^n.
\end{align}
Then by \eqref{ff1}, \eqref{ff2} and \eqref{eq:maximal B} with $0<r<\min\{p/r_w,q\}$ and $\theta$ from Lemma \ref{lem:FS} (b) we have
\begin{align*}
	\Vert S_\varphi f\Vert_{f^{p,q}_{\alpha,w}}
	&\lesssim \Big\Vert \Big(\sum_{j\ge 0} 2^{j\alpha q} \sum_{R\in\E_j} \big(\ip{f^{\varphi,R}}^* \Ind_R\big)^q\Big)^{1/q}\Big\Vert_{L^p_w} \\
	&\lesssim \Big\Vert \Big(\sum_{j\ge 0} 2^{j\alpha q} \Big(\sum_{P\in\E_{j+k}} \ip{f_{\varphi,P,k}}^* \Ind_P\Big)^q\Big)^{1/q}\Big\Vert_{L^p_w} \\
	&\lesssim  \Big\Vert \Big(\sum_{j\ge 0} \MM^\theta_r\Big(2^{j\alpha} \sum_{P\in\E_{j+k}} \ip{f_{\varphi,P,k}}\Ind_P\Big)^q\Big)^{1/q}\Big\Vert_{L^p_w}.
	\end{align*}
Now we invoke 	Lemma \ref{lem:FS} (b) and arrive at
\begin{align*}
\Vert S_\varphi f\Vert_{f^{p,q}_{\alpha,w}}
	\lesssim \Big\Vert \Big(\sum_{j\ge 0}\Big(2^{j\alpha} \sum_{P\in\E_{j+k}} \ip{f_{\varphi,P,k}}\Ind_P\Big)^q\Big)^{1/q}\Big\Vert_{L^p_w} 
	\le \Vert f\Vert_{F^{p,q}_{\alpha,w}},
\end{align*}
which completes the proof of part (b).

\medskip
\underline{\it The Besov case.}
As in the Triebel--Lizorkin case, part (c) follows from  Proposition \ref{prop: CRF} (see Remark \ref{rem:CRF}). 

\textbf{Part (a).} Assume that $F^{p,q}_{\alpha,w}$ has been defined using some admissible pair $(\wt{\vph}_\circ,\wt{\vph})$. 
Let $0<r<\min\{ p/r_w,1\}$ and $\theta$ be the number from Lemma \ref{lem:FS} (b). Then by Lemma~\ref{almost orthog} with $\sigma> \theta/r +\max\{n, n/r\}$, \eqref{eq:maximal A} with $a_R= |R|^{-1/2}|s_R|$, and Lemma \ref{lem:FS} (b) with $q=1$ we have
\begin{align*}
 \big\Vert \varphi_j(\sqrt{\LL})(T_{\psi} s)\big\Vert_{L^p_w}
	&\lesssim \sum_{k=j-2}^{j+2}\Big\Vert \MM^\theta_r\Big(\sum_{R\in\E_k} |R|^{-1/2}|s_R| \Ind_R\Big)\Big\Vert_{L^p_w}  \\
	&\lesssim \sum_{k=j-2}^{j+2} \Big\Vert \sum_{R\in\E_k} |R|^{-1/2}|s_R| \Ind_R\Big\Vert_{L^p_w} \\
	&\lesssim \sum_{k=j-2}^{j+2} \Big(\sum_{R\in\E_k} \big(w(R)^{1/p}|R|^{-1/2} |s_R| \big)^p\Big)^{1/p}.
\end{align*}
From the preceding estimate it follows readily that
\begin{align*}
	\Vert T_{\psi}s\Vert_{B^{p,q}_{\alpha,w}}
	\lesssim \Big\{ \sum_{j\ge 0} 2^{j\alpha q}\Big(   \sum_{k=j-2}^{j+2} \Big(\sum_{R\in\E_k} \big(w(R)^{1/p}|R|^{-1/2} |s_R| \big)^p\Big)^{1/p}\Big)^q \Big\}^{1/q}
	\lesi \Vert s\Vert_{b^{p,q}_{\alpha,w}}.
	\end{align*}
	
	\textbf{Part (b).} As in the Triebel--Lizorkin case assume that $F^{p,q}_{\alpha,w}$ has been defined using the admissible pair $(\vph_\circ,\vph)$ from the definition of $S_\vph$. Since $\varphi_j(\sqrt{\LL})f\in V_{4^j}$, we may  apply Lemma \ref{lem:plancherel-polya} to obtain
	\begin{align*}
		\Vert S_\varphi f\Vert_{b^{p,q}_{\alpha,w}}
		\lesssim \Big\{ \sum_{j\ge 0} 2^{j\alpha q}\Big(\sum_{R\in\E_j} w(R)\, |\varphi_j(\sqrt{\LL})f(x_R)|^p \Big)^{q/p}\Big\}^{1/q} 
		\lesi \Vert f\Vert_{B^{p,q}_{\alpha,w}}.
	\end{align*}

This concludes the Besov case and also completes the proof of Theorem \ref{thm:frames}.

\section{Embeddings for weighted Hermite spaces}\label{sec:embedding}
In this section we characterize continuous Sobolev-type embeddings  for weighted Besov and Triebel--Lizorkin spaces associated the Hermite operator. The main result in this section is the following.

\begin{Theorem}[Embedding for sequence spaces]\label{thm:embedding2}
Let $\gamma>0$ and $w$ be any weight. Then the following holds. 
\begin{enumerate}[\upshape(a)]
\item For any $\alpha_1\le \alpha_2$ and $0<p_1, p_2, q_1,q_2\le \infty$ with $0<q_2\le q_1\le \infty$ and $\alpha_1-\f{\gamma}{p_1} = \alpha_2 - \f{\gamma}{p_2}$, we have
$$ b^{p_2,q_2}_{\alpha_2,w}\embed b^{p_1,q_1}_{\alpha_1,w}$$
if and only if $w$ satisfies the lower bound property \eqref{eq:HLB} of order $\gamma$.

\item If in addition $w\in \AH_\infty$, then for any $\alpha_1\le \alpha_2$, $0<p_1, p_2<\infty$ and $0< q_1,q_2\le \infty$ with $\alpha_1-\f{\gamma}{p_1} = \alpha_2 - \f{\gamma}{p_2}$, we have
$$ f^{p_2,q_2}_{\alpha_2,w}\embed f^{p_1,q_1}_{\alpha_1,w}$$
if and only if $w$ satisfies the lower bound property \eqref{eq:HLB} of order $\gamma$.

\end{enumerate}
\end{Theorem}

Theorem \ref{thm:embedding2} will be used to prove Theorem \ref{thm:embedding1} in Section \ref{sec:embedding1proof}.
The proof of Theorem \ref{thm:embedding2} will be given in Section \ref{sec:embedding2proof}.
Before turning to its proof we state a crucial fact concerning the the lower bound property \eqref{eq:HLB}. Recall that the function $\cro(\cdot)$ was defined in \eqref{rho}.
\begin{Proposition}[Lower bound for tiles and balls]\label{prop: lb}
Suppose that $\gamma>0$ and let $w$ be any weight. The following are equivalent.
\begin{enumerate}[\upshape(a)]
	\item There exists $C>0$ such that for every $j\ge 0$ and $R\in\E_j$ we have 
	$$w(R) \ge C \,2^{-j\gamma}.$$
	\item There exists $\wt{C}>0$ such that for every ball $B$ with $0<r_B\le \cro(x_B)$ we have $$w(B)\ge \wt{C}\, r_B^\gamma.$$ 
\end{enumerate}
\end{Proposition}
\begin{proof}[Proof of Proposition \ref{prop: lb}]
	Proof of (a) $\Rightarrow$ (b).
	
Fix a ball $B$ with $0<r_B\le \cro(x_B)$. We observe that there exists $j\ge 0$ with 
\begin{align}\label{lb1}
	2^j\le |x_B|_\infty <2^{j+1}
\end{align}
and thus it follows that $r_B\le 2^{-j}$. Furthermore, there exists $k\ge 0$ with 
\begin{align}\label{lb2}
2^{-(j+k+1)}<r_B\le 2^{-(j+k)}.
\end{align}
Next choose an integer $\ell$ such that 
\begin{align}\label{lb3}
	2^\ell \ge \max\Big\{4 c_0 \sqrt{n}, \Big(\f{c_2}{4\delta_\star}\Big)^3\Big\}
\end{align}
Now consider tiles from $\E_{j+k+\ell}$. Firstly there is a tile $R$ from the collection $\E_{j+k+\ell}$ that contains the centre of $B$. Indeed from \eqref{lb1} and our assumption on $\ell$ we have $ |x_B|_\infty \le 2^{j+k+\ell}$ and hence, by property \eqref{Hermite tiles 3}, it follows that
$$ x_B\in Q\big(0,2^{j+k+\ell}\big) \subset \Q_{j+k+\ell}.$$
Secondly such a tile $R$ is contained in a cube of diameter approximately $2^{-(j+k+\ell)}$. To see this we note that by  \eqref{Hermite tiles 1-2} and  \eqref{lb3},
$$ |x_B-x_R|\le c_2 2^{-(j+k+\ell)/3}\le 4\delta_\star.$$
It follows that $|x_R|_\infty \le (1+4\delta_\star)2^{j+k+\ell+1}$ and so by \eqref{Hermite tiles 0}
\begin{align}\label{lb4} 
R \subset Q\big(x_R, c_0 2^{-(j+k+\ell)}\big). 
\end{align}
Thirdly $R$ is contained in $B$. Indeed since $x_B\in R$ and by \eqref{lb4}, \eqref{lb3} and \eqref{lb2} we have
\begin{align*}
	\text{diam}(R)=2c_0\sqrt{n}\,2^{-(j+k+\ell)}\le 2^{-(j+k+1)}< r_B.
\end{align*}
Hence $R\subset B$.
Combining the above three facts on $B$ and $R$ we conclude, by \eqref{lb2},
\begin{align*}
	w(B)\ge w(R) \ge C\, 2^{-(j+k+\ell)\gamma} \ge \wt{C}\, r_B^\gamma
\end{align*}
 with $\wt{C}=C\,2^{-\ell\gamma}$.

We turn to the proof of (b) $\Rightarrow$ (a).
Fix $j\ge 0$ and let $R$ be a tile from $\E_j$. Then consider the ball $B(x_R,r)$, where 
\begin{align*}
r:= \min\Big\{ c_1, \f{1}{1+c_3}\Big\}\,2^{-j}.
\end{align*}
Then firstly we have $B\subset R$ by \eqref{Hermite tiles 1-2}. Secondly we know that $|x_R|\le c_3 2^j$ by \eqref{Hermite tiles 3}. Thus $ r\le \cro(x_R)$ 
and so  $B(x_R, r)$ meets the hypothesis of Proposition \ref{prop: lb} (b). Thus we have
$$ w(R)\ge w(B) \ge \wt{C}\, r^\gamma = C\, 2^{-j\gamma},$$
where $C = \wt{C}\,\min\{ c_1, (1+c_3)^{-1}\}^\gamma$.
\end{proof}

We are now ready to give the proofs of our embedding theorems. 
\subsection{Proof of Theorem \ref{thm:embedding2}}\label{sec:embedding2proof}

\textbf{Part (a).} We first prove necessity of the lower bound property.
Suppose that for some $C>0$ we have
\begin{align}\label{bs1} \Vert s\Vert_{b^{p_1,q_1}_{\alpha_1,w}}\le C \Vert s\Vert_{b^{p_2,q_2}_{\alpha_2,w}}\end{align}
for all sequences . 
Fix $j\ge 0$ and a tile $R \in \E_j$ and consider the sequence with $s_R=1$ and $s_{\wt{R}}=0$ for any tiles $\wt{R}\ne R$. Then for any $p,q$ and $\alpha$ we have 
$$ \Vert s\Vert_{b^{p,q}_{\alpha, w}} =\Big\{ 2^{j\alpha q} \big(w(R)|R|^{-p/2}\big)^{q/p}\Big\}^{1/q}= 2^{j\alpha}w(R)^{1/p}|R|^{-1/2}.$$
Then \eqref{bs1} implies that 
$$ 2^{j\alpha_1}w(R)^{1/p_1} |R|^{-1/2} \lesssim 2^{j\alpha_2}w(R)^{1/p_2}|R|^{-1/2},$$
which in turn gives
$$ w(R)^{1/p_2-1/p_1} \gtrsim  2^{-j\gamma(1/p_2-1/p_1)}.$$
Now noting that $p_2\le p_1$ we obtain
$w(R)\gtrsim \, 2^{-j\gamma}$, and since this holds for any tile $R$, then Proposition \ref{prop: lb} yields the lower bound property for $w$ with order~$\gamma$. 

We now prove sufficiency for part (a). Since $q_1\ge q_2$ and $p_1\ge p_2$ we may apply the $q_2/q_1$-triangle inequality  followed by the $p_2/p_1$-triangle inequality to obtain
\begin{align*}
	\Vert s\Vert_{b^{p_1,q_1}_{\alpha_1,w}}
	&\le \bigg\{ \sum_{j\ge 0} 2^{j\alpha_1 q_2}\Big(\sum_{R\in\E_j} \big(w(R)^{1/p_1} |R|^{-1/2}|s_R|\big)^{p_1} \Big)^{q_2/p_1}\bigg\}^{1/q_2} \\
	&\le \bigg\{ \sum_{j\ge 0} 2^{j\alpha_1 q_2}\Big(\sum_{R\in\E_j} \big(w(R)^{1/p_1-1/p_2}w(R)^{1/p_2} |R|^{-1/2}|s_R|\big)^{p_2} \Big)^{q_2/p_2}\bigg\}^{1/q_2}.
	\end{align*}
Now the lower bound property of $w$ and Proposition~\ref{prop: lb} allow us to control the last expression by a constant multiple of
\begin{align*}	
		\bigg\{ \sum_{j\ge 0} 2^{j\alpha_1 q_2}\Big(\sum_{R\in\E_j} \big(2^{-j\gamma(1/p_1-1/p_2)}w(R)^{1/p_2} |R|^{-1/2}|s_R|\big)^{p_2} \Big)^{q_2/p_2}\bigg\}^{1/q_2},
\end{align*}	
which equals $\Vert s\Vert_{b^{p_2,q_2}_{\alpha_2,w}}$.

\textbf{Part (b).} Necessity can be done in the same way as part (a), after observing that for the same sequence $s$ defined in (a) we have
$$ \Vert s\Vert_{f^{p,q}_{\alpha, w}} = \Big(\int_R 2^{j\alpha p} |R|^{-p/2} w(x)\,dx\Big)^{1/p}=2^{j\alpha}w(R)^{1/p}|R|^{-1/2}$$
for any $p,q,\alpha$ and $w$.

\medskip
We turn to sufficiency of the lower bound property for part (b). 
We shall prove that there exists $C>0$ with 
\begin{align}\label{fs2}
\Vert s\Vert_{f^{p_1,q_1}_{\alpha_1,w}}\le C \Vert s\Vert_{f^{p_2,q_2}_{\alpha_2,w}}
\end{align}
for all sequences $s\in f^{p_2,q_2}_{\alpha_2,w}$.  Without loss of generality we may assume that $s=~\{s_R\}_R$ is a sequence with $\Vert s\Vert_{f^{p_2,q_2}_{\alpha_2,w}}=1$. 
To simplify notation we denote 
\begin{align*}
	F_j(\alpha, q, x):= \sum_{R\in\E_j} 2^{j\alpha q}\big(|R|^{-1/2} |s_R| \Ind_R(x)\big)^q.
\end{align*} 
We claim that \eqref{fs2} will follow from the following two estimates: there exists $\wt{\beta}>0$ such that for each $x\in\RR^n$ and any $N\ge 0$,
\begin{align}\label{fs3}
	\Big\{ \sum_{j=0}^N F_j(\alpha_1,q_1,x)\Big\}^{1/q_1} \le \wt{\beta}\,2^{N\gamma/p_1}; 
\end{align}
and  for any $N\ge -1$,
\begin{align}\label{fs4}
	\Big\{ \sum_{j=N+1}^\infty F_j(\alpha_1,q_1,x)\Big\}^{1/q_1} \le 2^{-N \gamma(\f{1}{p_2}-\f{1}{p_1})} \Big\{ \sum_{j=N+1}^\infty F_j(\alpha_2,q_2,x)\Big\}^{1/q_2}.
\end{align}
Let us show how \eqref{fs3} and \eqref{fs4} lead to \eqref{fs2}. We first discretize the left hand side of \eqref{fs2} using the well known representation of the $L^p$ norm as follows.
\begin{align*}
	\Vert s\Vert^{p_1}_{f^{p_1,q_1}_{\alpha_1,w}}
	&= p_1 \int_0^\infty t^{p_1-1} w\Big(\big\{x\in\RR^n: \Big(\sum_{j\ge 0} F_j(\alpha_1,q_1,x)\Big)^{1/q_1}>t\big\}\Big)\,dt =:I + II, 
	\end{align*}
where
\begin{align*}
	I &= p_1\int_0^\beta   t^{p_1-1} w\Big(\big\{x\in\RR^n: \Big(\sum_{j\ge 0} F_j(\alpha_1,q_1,x)\Big)^{1/q_1}>t\big\}\Big)\,dt,  \\
	II &=p_1 \sum_{N=0}^\infty \int_{\beta 2^{N\gamma/p_1}}^{\beta2^{(N+1)\gamma/p_1}}   t^{p_1-1} w\Big(\big\{x\in\RR^n: \Big(\sum_{j\ge 0} F_j(\alpha_1,q_1,x)\Big)^{1/q_1}>t\big\}\Big)\,dt, 
\end{align*}
and $\beta=2\wt{\beta} 2^{1/q_1}$ with $\wt{\beta}$ being the constant from \eqref{fs3}. 

We estimate $I$ by applying \eqref{fs4} with $N=-1$ along with the substitution $t = 2^{\gamma(1/p_2-1/p_1)}u$. Thus,
\begin{align*}
	I &\le p_1\int_0^\beta   t^{p_1-1} w\Big(\big\{x\in\RR^n: 2^{\gamma(\f{1}{p_2}-\f{1}{p_1})} \Big(\sum_{j\ge 0} F_j(\alpha_2,q_2,x)\Big)^{1/q_2}>t\big\}\Big)\,dt \\
	&=p_1\int_0^\beta   t^{p_1-1} w\Big(\big\{x\in\RR^n: \Big(\sum_{j\ge 0} F_j(\alpha_2,q_2,x)\Big)^{1/q_2}>2^{-\gamma(\f{1}{p_2}-\f{1}{p_1})} t\big\}\Big)\,dt \\
	&=p_12^{\gamma(\f{p_1}{p_2}-1)}\int_0^{\beta 2^{-\gamma(\f{1}{p_2}-\f{1}{p_1})}}   u^{p_1-1} w\Big(\big\{x\in\RR^n: \Big(\sum_{j\ge 0} F_j(\alpha_2,q_2,x)\Big)^{1/q_2}> u\big\}\Big)\,du.
\end{align*}
Now from inequality $u^{p_1}\le u^{p_2}\big(\beta 2^{-\gamma(1/p_2-1/p_1)}\big)^{p_1-p_2}$, we have
\begin{align*}
	I&\lesssim p_2\int_0^{\beta 2^{-\gamma(\f{1}{p_2}-\f{1}{p_1})}}   u^{p_2-1} w\Big(\big\{x\in\RR^n: \Big(\sum_{j\ge 0} F_j(\alpha_2,q_2,x)\Big)^{1/q_2}> u\big\}\Big)\,du
	\lesssim \Vert s\Vert^{p_2}_{f^{p_2,q_2}_{\alpha_2,w}} = 1.
\end{align*}

For term $II$ we first employ \eqref{fs3} and \eqref{fs4} respectively to obtain
\begin{align*}
	II & \le p_1 \sum_{N=0}^\infty \int_{\beta 2^{N\gamma/p_1}}^{\beta2^{(N+1)\gamma/p_1}}   t^{p_1-1} w\Big(\big\{x\in\RR^n: \Big(\sum_{j\ge 0} F_j(\alpha_1,q_1,x)\Big)^{1/q_1}> 2^{-1/q_1}t\big\}\Big)\,dt \\
	& \le p_1 \sum_{N=0}^\infty \int_{\beta 2^{N\gamma/p_1}}^{\beta2^{(N+1)\gamma/p_1}}   t^{p_1-1} w\Big(\big\{x\in\RR^n: \Big(\sum_{j=N+1}^{\infty} F_j(\alpha_1,q_1,x)\Big)^{1/q_1}> \tfrac{2^{-1/q_1}}{2}t\big\}\Big)\,dt \\
	& \le p_1 \sum_{N=0}^\infty \int_{\beta 2^{N\gamma/p_1}}^{\beta2^{(N+1)\gamma/p_1}}   t^{p_1-1} w\Big(\big\{x\in\RR^n: \Big(\sum_{j=N+1}^{\infty} F_j(\alpha_2,q_2,x)\Big)^{1/q_2}> 2^{N\gamma(\f{1}{p_2}-\f{1}{p_1})}\tfrac{2^{-1/q_1}}{2}t\big\}\Big)\,dt.
\end{align*}
Now note that $t\sim \beta 2^{N\gamma/p_1}$ implies $2^{N\gamma(\f{1}{p_2}-\f{1}{p_1})}\tfrac{2^{-1/q_1}}{2}t\sim \beta' t^{p_1/p_2}$ with $\beta' = \beta^{1-p_1/p_2}2^{-\gamma/p_2}$. This estimate, along with the substitution $u=\beta' t^{p_1/p_2}$, gives 
\begin{align*}
II&\le p_1 \sum_{N=0}^\infty \int_{\beta 2^{N\gamma/p_1}}^{\beta2^{(N+1)\gamma/p_1}}   t^{p_1-1} w\Big(\big\{x\in\RR^n: \Big(\sum_{j=N+1}^{\infty} F_j(\alpha_2,q_2,x)\Big)^{1/q_2}>  \beta' t^{p_1/p_2}\big\}\Big)\,dt \\
	&\lesssim p_2\int_0^\infty   u^{p_2-1} w\Big(\big\{x\in\RR^n: \Big(\sum_{j\ge 0}F_j(\alpha_2,q_2,x)\Big)^{1/q_2}>  u\big\}\Big)\,du,
\end{align*}
which equals $\Vert s\Vert^{p_2}_{f^{p_2,q_2}_{\alpha_2,w}}$.

It remains to prove estimates \eqref{fs3} and \eqref{fs4}. We begin with \eqref{fs3}. 
The idea is to apply almost orthogonality to estimate $|s_R|$. Let $\{\varphi_j\}_{j\ge 0}$ and $\{\psi_j\}_{j\ge 0}$ be admissible systems so that \eqref{eq:CRF1} holds. Then we have 
$$ 
s_R = \ip{T_\psi s, \varphi_R} = \tau_R^{1/2}\varphi_j(\sqrt{\LL})(T_\psi s)(x_R),
$$
and so by Lemma \ref{almost orthog} we have for any $\sigma \ge 1$, 
\begin{align}\label{fs5}
	|s_R|
	\le |\tau_R^{1/2}|  \big|\varphi_j(\sqrt{\LL})(T_\psi s)(x_R)\big|
	\lesssim |R|^{1/2}\sum_{k=j-2}^{j+2}\sum_{P\in\E_k} \f{|P|^{-1/2}|s_P|}{(1+2^k|x_R-x_P|)^\sigma}.
\end{align}
Now fix an $r\in (0, \min\{p_2/r_w, q_2\})$ and let $\theta$ be the number in Lemma \ref{lem:FS}. Then applying \eqref{eq:maximal A} to \eqref{fs5} with $\sigma>\theta/r + \max\{n,n/r\}$ and $a_P=|P|^{-1/2}|s_P|$, we have
\begin{align}\label{fs7}
	|R|^{-1/2}|s_R|
	\lesssim \sum_{k=j-2}^{j+2}\sum_{P\in\E_k} \f{|P|^{-1/2}|s_P|}{(1+2^k|x_R-x_P|)^\sigma}
	\lesssim \sum_{k=j-2}^{j+2}\MM^\theta_r\Big(\sum_{P\in\E_k} |P|^{-1/2}|s_P| \Ind_P\Big)(x')
\end{align}
for any   $x'\in Q(x,2c_4 2^{-k})$ and $x\in R$.

Write $\wt{Q}:=Q(x,2c_4 2^{-k})$. Then for each $k\in\{j-2,j-1,j,j+1,j+2\}$ we have, by Lemma~\ref{lem:FS} (b),
\begin{align*}
\inf_{x'\in \wt{Q}} \MM^\theta_r\big(F_k(\alpha_2,1,\cdot)\big)(x')
	&\le\Big\{ \inf_{x'\in \wt{Q}} \Big[ \sum_{k\ge 0}\MM^\theta_r\big(F_k(\alpha_2,1,\cdot)\big)(x')^{q_2}\Big]^{\f{p_2}{q_2}}\Big\}^{\f{1}{p_2}}\\  
	&\le \Big\{\f{1}{w(\wt{Q})}\int_{\wt{Q}} \Big( \sum_{k\ge 0} \MM^\theta_r\big(F_k(\alpha_2,1,\cdot)\big)(y)^{q_2}\Big)^{\f{p_2}{q_2}}w(y)\,dy\Big\}^{\f{1}{p_2}}\\
	&\le w(\wt{Q})^{-1/p_2} \Big\Vert \Big(\sum_{k\ge 0} \Big| \MM^\theta_r\big(F_k(\alpha_2,1,\cdot)\big)\Big|^{q_2} \Big)^{\f{r}{q_2}}\Big\Vert^{\f{1}{r}}_{L^{\f{p_2}{r}}_w} \\
	&\lesi w(\wt{Q})^{-1/p_2} \Big\Vert \Big(\sum_{k\ge 0} \Big| \big(F_k(\alpha_2,1,\cdot)\big)\Big|^{q_2} \Big)^{\f{r}{q_2}}\Big\Vert^{\f{1}{r}}_{L^{\f{p_2}{r}}_w}.
\end{align*}
Thus noting that
\begin{align*}
\Big\Vert \Big(\sum_{k\ge 0} \Big| \big(F_k(\alpha_2,1,\cdot)\big)\Big|^{q_2} \Big)^{\f{r}{q_2}}\Big\Vert^{\f{1}{r}}_{L^{\f{p_2}{r}}_w}
=\Big\Vert \Big(\sum_{k\ge 0}F_k(\alpha_2,q_2,\cdot)\Big)^{\f{1}{q_2}}\Big\Vert_{L^{p_2}_w}
	= \Vert s\Vert_{f^{p_2,q_2}_{\alpha_2,w}} 
	=1,
\end{align*}
we arrive at the estimate
\begin{align}\label{fs8}
\inf_{x'\in \wt{Q}} \MM^\theta_r\Big(\sum_{P\in\E_k} |P|^{-1/2}|S_P| \Ind_P\Big)(x') 
=2^{-k\alpha_2}\inf_{x'\in \wt{Q}} \MM^\theta_r\big(F_k(\alpha_2,1,\cdot)\big)(x')
	\lesssim 2^{-k\alpha_2} w(\wt{Q})^{-1/p_2}.
\end{align}

Therefore by combining \eqref{fs7} and \eqref{fs8}, and observing that because $\sum_{R\in\E_j}\Ind_R =\Ind_{\Q_j}\le 1$ and $2^{-k\alpha_2}\sim 2^{-j\alpha_2}$, we have for each $x\in\RR^n$,
\begin{align}\label{fs9}
	&\Big\{ \sum_{j=0}^N F_j(\alpha_1,q_1,x)\Big\}^{1/q_1} \notag\\  \notag
	&\qquad=\Big\{ \sum_{j=0}^N \sum_{R\in\E_j}2^{j\alpha_1 q_1} \big(|R|^{-1/2}|s_R|\Ind_R(x)\big)^{q_1} \Big\}^{1/q_1} \\ \notag
	&\qquad\lesssim \Big\{ \sum_{j=0}^N \sum_{R\in\E_j} 2^{j\alpha_1 q_1} \Ind_R(x) \bigg(\sum_{k=j-2}^{j+2}\inf_{x'\in\wt{Q}}\MM^\theta_r\Big(\sum_{P\in\E_k} |P|^{-1/2}|s_P| \Ind_P\Big)(x')\bigg)^{q_1} \Big\}^{1/q_1} \\  \notag
	&\qquad\lesssim \Big\{ \sum_{j=0}^N \sum_{R\in\E_j} 2^{j\alpha_1 q_1} \Ind_R(x) \bigg(\sum_{k=j-2}^{j+2} 2^{-k\alpha_2} w(\wt{Q})^{-1/p_2}\bigg)^{q_1} \Big\}^{1/q_1} \\ 
	&\qquad\lesssim \Big\{ \sum_{j=0}^N 2^{jq_1(\alpha_1-\alpha_2)}\bigg(\sum_{k=j-2}^{j+2}  w(\wt{Q})^{-1/p_2}\bigg)^{q_1} \Big\}^{1/q_1}.
\end{align}
To complete the estimate we apply the lower bound condition to handle $w(\wt{Q})^{-1/p_2}$. First note that if $x\notin Q(0,c_32^j)$ then $\Ind_R(x)=0$ for every $R\in\E_j$ and so estimate \eqref{fs3} holds trivially. Thus we may assume $x\in Q(0,c_32^j)$. Let us take the ball $B$ with $x_B=x$ and 
$$ r_B = \min\Big\{\f{1}{1+c_3}, c_4\Big\}\,2^{-j}.$$
Then we have $B\subset \wt{Q}$ because $|k-j|\le 2$ implies $r_B\le 2c_42^{-k}$. We also have $|x|_\infty \le c_32^j$ because $x\in Q(0,c_32^j)$, and so $r_B\le \cro(x_B)$. Thus we may apply the lower bound condition to $B$ and obtain
$$ w(\wt{Q})\ge w(B)\gtrsim r_B^\gamma \sim 2^{-j\gamma}.$$
Inserting this into the final line of \eqref{fs9}, we have for each $x\in Q(0,c_32^j)$,
\begin{align*}
	\Big\{ \sum_{j=0}^N F_j(\alpha_1,q_1,x)\Big\}^{1/q_1} 
	\lesssim \Big\{ \sum_{j=0}^N 2^{jq_1(\alpha_1-\alpha_2)} 2^{j\gamma q_1/p_2} \Big\}^{1/q_1} 
	= \Big\{\sum_{j=0}^N 2^{j\gamma q_1/p_1}\Big\}^{1/q_1} 
	\le \wt{\beta}\, 2^{N \gamma/p_1}.
\end{align*}
This proves \eqref{fs3}, after taking into account that that the estimate holds trivially if $x\notin Q(0,c_32^j)$. 

We turn to the proof of \eqref{fs4}. Now since all the tiles in $\E_j$ are disjoint we have
\begin{align}\label{fs10}
	\Big\{ \sum_{j=N+1}^\infty F_j(\alpha_1,q_1,x)\Big\}^{1/q_1}
	=
	\Big\{\sum_{j=N+1}^\infty 2^{jq_1(\alpha_1-\alpha_2)} F_j(\alpha_2,1,x)^{q_1}\Big\}^{1/q_1}.
\end{align}
If $q_1\ge q_2$ we apply the $q_2/q_1$-triangle inequality, but if $q_1\le q_2$ then we apply H\"older's inequality with exponent $q_2/q_1$ to the the terms $2^{jq_2(\alpha_1-\alpha_2)/2} F_j(\alpha_2,1,x)$ and $2^{jq_2(\alpha_1-\alpha_2)/2}$. In either case \eqref{fs10} implies
\begin{align*}
	\Big\{ \sum_{j=N+1}^\infty F_j(\alpha_1,q_1,x)\Big\}^{1/q_1}
	&\le \Big\{\sum_{j=N+1}^\infty 2^{jq_2(\alpha_1-\alpha_2)} F_j(\alpha_2,1,x)^{q_2}\Big\}^{1/q_2} \\
	&\lesssim 2^{N(\alpha_1-\alpha_2)} \Big\{ \sum_{j=N+1}^\infty F_j(\alpha_2,q_2,x)\Big\}^{1/q_2},
\end{align*}
which yields \eqref{fs4}.

This completes the proof of \eqref{fs3} and \eqref{fs4} and hence of Theorem \ref{thm:embedding2}.

\subsection{Proof of Theorem \ref{thm:embedding1}}\label{sec:embedding1proof}
We can now bring together the material developed throughout the paper to prove our main result stated in the Introduction.
The theorem follows by combining our frame decompositions with our embedding result for the sequence spaces. 

Suppose firstly that the lower bound property \eqref{eq:HLB} holds. Then  Theorem \ref{thm:embedding2} in conjunction with Theorem \ref{thm:frames} gives
\begin{align*}
		 \Vert f \Vert_{A^{p_1,q_1}_{\alpha_1,w}(\LL)} 
		 \sim  \Vert \{\ip{f,\varphi_R}\}_R\Vert_{a^{p_1,q_1}_{\alpha_1,w}(\LL)}  
		 \lesssim \Vert \{\ip{f,\varphi_R}\}_R\Vert_{a^{p_2,q_2}_{\alpha_2,w}(\LL)}  
		 \sim \Vert f \Vert_{A^{p_2,q_2}_{\alpha_2,w}(\LL)},
 \end{align*}
 which shows that $A^{p_2,q_2}_{\alpha_2,w}\embed A^{p_1,q_1}_{\alpha_1,w}$.
 
 On the other hand, if \eqref{eq:Femb} or \eqref{eq:Bemb} holds then  invoking Theorem \ref{thm:frames}  again we obtain
 \begin{align*}
  \Vert s\Vert_{a^{p_1,q_1}_{\alpha_1,w}(\LL)} 
  \sim  \Vert T_\psi s \Vert_{A^{p_1,q_1}_{\alpha_1,w}(\LL)} 
    \lesssim  \Vert T_\psi s \Vert_{A^{p_2,q_2}_{\alpha_2,w}(\LL)} 
    \sim   \Vert s\Vert_{a^{p_2,q_2}_{\alpha_2,w}(\LL)}.
  \end{align*}
 Thus by Theorem \ref{thm:embedding2} we see that \eqref{eq:HLB} holds, concluding the proof of Theorem \ref{thm:embedding1}.

\bigskip

\textbf{Acknowledgment.} The first named author was partly supported by Ho Chi Minh City University of Education under the B2020 project entitled ``Singular integrals with non-smooth kernels on function spaces associated with differential operators''. The second named author was supported by ARC grant DP170101060.

\end{document}